\documentclass[12pt]{article}

\usepackage[latin1]{inputenc}

\usepackage{pgf,tikz}
\usetikzlibrary{arrows}
\usetikzlibrary{patterns}

\usepackage{subfigure}

\usepackage{amssymb}

\usepackage{amsmath,epsfig}
\def\WE{\textrm{WE}}
\def\NS{\textrm{NS}}
\def\NE{\textrm{NE}}
\def\SE{\textrm{SE}}
\def\NW{\textrm{NW}}
\def\SW{\textrm{SW}}
\def\XY{\textrm{XY}}
\def\X{\textrm{X}}
\def\Y{\textrm{Y}}

\newtheorem{Theorem}{Theorem}
\newtheorem{Corollary}[Theorem]{Corollary}

\newtheorem{Lemma}[Theorem]{Lemma}

\newcommand{\C}{{\mathbb C}}

\newcommand{\beq}{\begin{equation}}
\newcommand{\eeq}{\end{equation}}

\def\qed{\hfill$\Box$}

\usepackage{multirow}

\def\mybox{\hbox to 12.0pt}

\def\mybigbox{\hbox to 35.0pt}
\def\myverybigbox{\hbox to 60.0pt}

\def\ds{\displaystyle}    
\def\ov{\overline}

\def\wgt{{\rm wgt}}

\def\x{{\bf x}}
\def\y{{\bf y}}
\def\z{{\bf z}}

\def\a{{\bf a}}

\def\1{{\bf 1}}

\def\n{{\bf n}}
\def\sh{S} 
\def\T{{t}}  
\def\S{{s}}  
\def\P{{p}}  
\def\Q{{q}}
\def\+{\!\!+\!\!}
\def\red{\textcolor{red}}
\def\blue{\textcolor{blue}}
\def\magenta{\textcolor{magenta}}

\newcommand{\YT}[3]{
\vcenter{\hbox{
\begin{tikzpicture}[x={(0in,-#1)},y={(#1,0in)}] 
\foreach \k [count=\i] in {#3} {
 \foreach \e [count=\j] in \k {
  \draw (\i,\j) rectangle +(-1,-1);
  \draw (\i-0.5,\j-0.5) node {$#2\e$};
 }
}
\end{tikzpicture}
}}
}

\newcommand{\wideYT}[4]{
\vcenter{\hbox{
\begin{tikzpicture}[x={(0in,-#1)},y={(#2,0in)}] 
\foreach \k [count=\i] in {#4} {
 \foreach \e [count=\j] in \k {
  \draw (\i,\j) rectangle +(-1,-1);
  \draw (\i-0.5,\j-0.5) node {$#3\e$};
 }
}
\end{tikzpicture}
}}
}

\newcommand{\SYT}[3]{
\vcenter{\hbox{
\begin{tikzpicture}[x={(0in,-#1)},y={(#1,0in)}] 
\foreach \k [count=\i] in {#3} {
 \foreach \e [count=\j] in \k {
  \draw (\i,\j+\i-1) rectangle +(-1,-1);
  \draw (\i-0.5,\j+\i-1-0.5) node {$#2\e$};
 }
}
\end{tikzpicture}
}}
}

\newcommand{\wideSYT}[4]{
\vcenter{\hbox{
\begin{tikzpicture}[x={(0in,-#1)},y={(#2,0in)}] 
\foreach \k [count=\i] in {#4} {
 \foreach \e [count=\j] in \k {
  \draw (\i,\j+\i-1) rectangle +(-1,-1);
  \draw (\i-0.5,\j+\i-1-0.5) node {$#3\e$};
 }
}
\end{tikzpicture}
}}
}

\title{Tokuyama's Identity for Factorial Schur Functions}

\author{
Ang\`ele M. Hamel\thanks{ 
Department of Physics and Computer Science,
Wilfrid Laurier University,
Waterloo, Ontario, N2L 3C5, Canada ({\tt ahamel@wlu.ca})}
\and 
Ronald C. King\thanks{
Mathematical Sciences, University of Southampton, 
Southampton SO17 1BJ, England ({\tt R.C.King@soton.ac.uk})}}

\begin{document}

\maketitle

\begin{abstract}
A recent paper of Bump, McNamara and Nakasuji introduced a factorial version of Tokuyama's identity, expressing the partition function of  six vertex model as the product of a $t$-deformed Vandermonde and a Schur function.  Here we provide an extension of their result by exploiting the language of primed shifted tableaux, with its proof based on the use of non-interesecting lattice paths.
\end{abstract}

\section{Introduction}
\label{sec:intro}

Tokuyama's identity~\cite{Tokuyama}, which expresses a weighted sum over strict Gelfand-Tsetlin patterns~\cite{GT} as the product 
of a $t$-deformed Vandermonde determinant and a Schur function, was originally established for $GL(n,\C)$ and its 
associated root system of type $A_{n-1}$, 
but subsequently other Tokuyama-like identities have been derived for other groups and their root systems~\cite{Brubaker1,Brubaker3,HKWeyl}. 
One of the recent additions to this literature is the paper of Bump, McNamara and Nakasuji~\cite{BMN}, who 
extended the original Tokuyama identity in a way that expresses the partition function of the six vertex model as the product of a factorial Schur function and the same $t$-deformed Vandermonde as before by using a six-vertex model interpretation due to Lascoux~\cite{Lascoux} and McNamara~\cite{McNamara} and the repeated application of the Yang--Baxter 
equation~\cite{Brubaker1}. 

Here we provide a further generalisation involving more than just a single deformation parameter $t$. To this end 
we make use of the fact that both the original Tokuyama identity and that of Bump {\em et al.} can be expressed in a natural manner in terms of certain primed shifted tableaux. Weighting these tableaux by means of two sets of indeterminates
$\x=(x_1,x_2,\ldots,x_n)$ and $\y=(y_1,y_2,\ldots,y_n)$, together with a sequence of shift parameters $\a=(a_1,a_2\ldots)$,
enables us to establish the required generalisation, with a proof provided by means of a non-intersecting lattice path argument.

Tokuyama's identity can be expressed, with a slight change of notation, in the form:
\begin{equation}\label{eqn-Tok}
          \sum_{G\in{\cal G}^\lambda}\ \wgt(G)  =  \prod_{1\leq i<j\leq n} (x_i + tx_j )\  \  s_\mu(\x) \,,
\end{equation}
where $\lambda=\mu+\rho$, with $\mu$ a partition with no more than $n$ parts and $\rho=(n-1,n-2,\ldots,1,0)$.
Here $\x=(x_1,x_2,\ldots,x_n)$ and $t$ are independent parameters. On the left,
the sum is over all strict Gelfand-Tsetlin patterns $G$ whose top row is the strict partition $\lambda$
and $\wgt(G)$, which will be specified later.
The reader will recognize $\prod_{1\leq i<j\leq n}(x_i + tx_j)$ as the expansion of a Vandermonde determinant deformed by the parameter $t$. The term $s_\mu(\x)$ is a Schur function, defined for example in the texts by Littlewood~\cite{Littlewood} and by Macdonald~\cite{MacText}.
Tokuyama's identity can be considered to be a deformation of Weyl's character formula for the reductive Lie algebra $gl(n)$ of the general linear group $GL(n)$ since at $t=-1$ one can recover the expression for the irreducible character $s_\mu(\x)$ as the ratio of two alternants.

The theorem of Bump, McNamara, and Nakasuji~\cite{BMN} states, again with a slight change of notation, that 
\begin{equation}\label{eqn-BMN}
Z(\mathfrak{S}_{\lambda,t}^{\Gamma}) = \prod_{i<j}(tx_i + x_j) \ \ s_{\lambda} (\x|\a),
\end{equation}
where $s_{\lambda} (\x|\a)$ is a factorial Schur function defined in Section~\ref{sec:main}.
The first such factorial Schur function was defined by Biedenharn and Louck~\cite{Biedenharn} 
in terms of Gelfand-Tsetlin patterns in a slightly more restricted form
(see also Chen and Louck~\cite{Chen}), but given its more general form by
Goulden and Greene~\cite{GC} and Macdonald~\cite{Macdonald}, 
expressed this time in terms of column-strict, that is to say semistandard, tableaux,
with Macdonald also giving an alternative definition as a ratio of alternants. 
The term $Z(\mathfrak{S}_{\lambda,t}^{\Gamma})$ is the partition function of the six vertex model 
$\mathfrak{S}_{\lambda,t}^{\Gamma}$ with a particular choice of Boltzmann weights that will also be specified later in Section~\ref{sec:corollaries}.

The combinatorial identities (\ref{eqn-Tok}) and (\ref{eqn-BMN}) due to Tokuyama~\cite{Tokuyama} and 
Bump {\em et al.}~\cite{BMN} that we are trying to generalise here were stated in terms of
strict Gelfand-Tsetlin patterns and the partition function of the square ice six vertex model. 
That one is a generalisation of the other comes about through the bijective correspondence between 
these two sets of combinatorial objects, together with the fact that a factorial Schur function is 
a generalisation of a Schur function. 
Here we will show that a natural combinatorial setting for
both these identities is that of primed shifted tableaux and associated non-intersecting lattice paths.

The paper is organized as follows: Section~\ref{sec:background} provides background information on tableaux and primed shifted tableaux, including definitions; Section~\ref{sec:main} gives our main result along with its proof based on a sequence of lemmas that are proved in Section~\ref{sec:lem-proofs}; and finally in Section~\ref{sec:corollaries} a number of corollaries are derived, special cases of which are shown to include both Tokuyama's identity and that of Bump et al.

\section{Tableaux and primed tableaux}
\label{sec:background}


To proceed we introduce some notation regarding partitions and tableaux. 
For any positive integer $n$ the sequence $\lambda = (\lambda_1, \lambda_2, \ldots, \lambda_n)$ with $ \lambda_1 \geq \lambda_2 \geq \ldots \geq \lambda_n$ is a partition if each part $\lambda_i$ is a non-negative integer. Its length $\ell(\lambda)$ is the number of non-zero parts and its weight $|\lambda|$ is the sum of its parts. We say that the partition $\lambda$ is strict if the above inequalities are all strict, i.e. all the parts of $\lambda$ are distinct.

A partition $\lambda$ of length $\ell(\lambda)\leq n$ defines a Young diagram $F^\lambda$ consisting of an array of $|\lambda|$ boxes $(i,j)$ arranged in rows of lengths $\lambda_i$ for $i=1,2,\ldots,\ell(\lambda)$ with $j=1,2,\ldots,\lambda_i$. Adopting 
the (English) convention whereby $(i,j)$ are matrix coordinates, the rows of $F^\lambda$ are left-adjusted to a vertical line.
If $\lambda$ is strict then it also defines a shifted Young diagram $SF^\lambda$ in which the rows of $F^\lambda$ are shifted to the right and left-adjusted to a diagonal line with boxes $(i,j)$ at $j=i,i+1,\ldots,i+\lambda_i-1$ for 
$i=1,2,\ldots,\ell(\lambda)$.  
Both $F^\lambda$ and $SF^\lambda$ consist of columns top-adjusted to a horizontal line.

For example, we have 
\begin{equation}\label{eqn-F-SF}
F^{3221}\ = \ 
\YT{0.2in}{}{
 {{},{},{}},
 {{},{}},
 {{},{}},
 {{}}
}
\qquad\qquad\qquad
SF^{6431}\ = \
\SYT{0.2in}{}{
 {{},{},{},{},{},{}},
 {{},{},{},{}},
 {{},{},{}},
 {{}}
}
\end{equation}

Using these conventions we define three different kinds of tableaux: semistandard tableaux,
shifted tableaux and primed shifted tableaux~\cite{StemShifted}. We restrict our attention to partitions $\lambda$
of length $\ell(\lambda)\leq n$ and strict partitions $\lambda$ of length $\ell(\lambda)=n$
and work with alphabets $[\n]=\{1<2<\cdots<n\}$, $[\n']=\{1'<2'<\cdots<n'\}$ and $[\n,\n']=\{1'<1<2'<2<\cdots<n'<n\}$. 


First, for each partition $\lambda$ let ${\cal T}^\lambda[\n]$ be the set of all semistandard tableaux $T$ of shape $\lambda$ 
that are obtained by filling each box $(i,j)$ of $F^\lambda$ with an entry $\T_{ij}\in[\n]$ in all possible ways 
such that:
\begin{itemize}
\item[T1] entries weakly increase from left to right across rows;
\item[T2] entries strictly increase from top to bottom down columns.
\end{itemize}

Then, for each strict partition $\lambda$ let ${\cal S}^\lambda[\n]$ be the set of all shifted tableaux $S$ of shape $\lambda$ that are obtained by filling each box $(i,j)$ of $SF^\lambda$ with an entry $\S_{ij}\in[\n]$ in all possible ways such that:
\begin{itemize}
\item[S1] entries weakly increase from left to right across rows;
\item[S2] entries weakly increase from top to bottom down columns;
\item[S3] entries strictly increase down each diagonal from top-left to bottom-right.
\end{itemize}

Finally, for each strict partition $\lambda$ let ${\cal Q}^\lambda[\n,\n']$ be the set of all primed shifted tableaux $P$ of shape $\lambda$ that are obtained by filling each box $(i,j)$ of $SF^\lambda$ with 
an entry $p_{ij}\in[\n,\n']$ in all possible ways such that:
\begin{itemize}
\item[P1] entries weakly increase from left to right across rows;
\item[P2] entries weakly increase from top to bottom down columns;
\item[P3] at most one entry $k'$ appears in any row for each $k\geq1$;
\item[P4] at most one entry $k$ appear in any column for each $k\geq1$,
\end{itemize}
and let ${\cal P}^\lambda[\n,\n']$ be the subset of ${\cal Q}^\lambda[\n,\n']$ such that:
\begin{itemize}
\item[P5] no primed entries appear on the main diagonal.
\end{itemize}

For example, we have 
\begin{equation}\label{eqn-TSP}
T\ = \ 
\YT{0.2in}{}{
 {1,2,4},
 {2,3},
 {4,4},
 {5}
}
\qquad
S\ = \
\SYT{0.2in}{}{
 {1,1,2,2,3,4},
   {2,3,3,3},
     {3,4,4},
       {4}
}
\qquad
P\ = \
\SYT{0.2in}{}{
 {1,1,2',2,3',4},
   {2,3',3,3},
     {3,4',4},
       {4}
}
\end{equation}
with $T\in{\cal T}^{(3,2,2,1,0)}([\bf5])$, $S\in{\cal S}^{(6,4,3,1)}([\bf4])$ and $P\in{\cal P}^{(6,4,3,1)}([\bf4,\bf4'])$. 

\section{Main Result}
\label{sec:main}

Let $\x=(x_1,x_2,\ldots,x_n)$, $\y=(y_1,y_2,\ldots,y_n)$ and $\a=(a_0,a_1,a_2,\ldots)$ be sequences of independent parameters. 
Then each partition $\lambda$ specifies not only the Schur function~\cite{Littlewood,MacText} 
\begin{equation}\label{eqn-Schur}
   s_\lambda(\x) = \sum_{T\in{\cal T}^\lambda(\n)}\ \prod_{\overset{(i,j)\in F^\lambda}{\T_{ij}\in[\n]}}\  x_{\T_{ij}}
\end{equation}
but also the factorial Schur function~\cite{GC,Macdonald}
\begin{equation}\label{eqn-factSchur}
   s_\lambda(\x|\a) = \sum_{T\in{\cal T}^\lambda(\n)}\ \prod_{\overset{(i,j)\in F^\lambda}{\T_{ij}\in[\n]}}\  (x_{\T_{ij}}+a_{\T_{ij}+j-i})\,.
\end{equation}	
Similarly, each strict partition $\lambda$ specifies not only the generalised Schur $P$ and $Q$-functions~\cite{HKbij} 
\begin{align}
   P_\lambda(\x;\y) &=\ds \sum_{P\in{\cal P}^\lambda(\n,\n')}
	                 \ \prod_{\overset{(i,j)\in F^\lambda}{\P_{ij}\in[\n]}}\  x_{\P_{ij}}
					         \ \prod_{\overset{(i,j)\in F^\lambda}{\P_{ij}\in[\n']}} \ y_{|\P_{ij}|}\,;\hfill\label{eqn-PSchur} \\ 
   Q_\lambda(\x;\y) &=\ds \sum_{P\in{\cal Q}^\lambda(\n,\n')}
	                 \ \prod_{\overset{(i,j)\in F^\lambda}{\P_{ij}\in[\n]}}\  x_{\P_{ij}}
					         \ \prod_{\overset{(i,j)\in F^\lambda}{\P_{ij}\in[\n']}} \ y_{|\P_{ij}|}\,,\hfill\label{eqn-QSchur} 
\end{align}	
but also the factorial generalised Schur $P$ and $Q$-functions introduced here for the first time in the form   
\begin{align}
   P_\lambda(\x;\y|\a) &=\ds  \sum_{P\in{\cal P}^\lambda(\n,\n')}
	                 \ \prod_{\overset{(i,j)\in F^\lambda}{\P_{ij}\in[\n]}}\  (x_{\P_{ij}}+a_{j-i})
	                 \ \prod_{\overset{(i,j)\in F^\lambda}{\P_{ij}\in[\n']}} \ (y_{|\P_{ij}|}-a_{j-i})~~\hbox{with $a_0=0$}; \hfill \label{eqn-factPSchur}\\ 
   Q_\lambda(\x;\y|\a) &=\ds  \sum_{P\in{\cal Q}^\lambda(\n,\n')}
	                 \ \prod_{\overset{(i,j)\in F^\lambda}{\P_{ij}\in[\n]}}\  (x_{\P_{ij}}+a_{j-i})
	                 \ \prod_{\overset{(i,j)\in F^\lambda}{\P_{ij}\in[\n']}} \ (y_{|\P_{ij}|}-a_{j-i})\,,\hfill\label{eqn-factQSchur} 
\end{align}
where in both cases $|\P_{ij}|=k$ if $\P_{ij}=k'$. It is notable here that the index on each $a$ is independent of $\P_{ij}$,
unlike the factorial Schur function case.

Note also that if we set $\y=\x$ and replace $\a$ by $-\a$ in $P_\lambda(\x;\y|\a)$ and $Q_\lambda(\x;\y|\a)$ they reduce 
to the factorial functions $P_\lambda(\x\,|\,0,\a)$ 
and $Q_\lambda(\x\,|\,\a)$ of Ikeda {\em et al.}, see section 4.2 of~\cite{ikeda}, that have been shown to be expressible combinatorially
in terms of 
primed shifted tableaux by Ivanov, as exemplified in his Theorem 2.11~\cite{ivanov}.
\bigskip

We may now state our main result:


\begin{Theorem}\label{thm-main}
Let $\mu$ be a partition of length $\ell(\mu)\leq n$ and $\delta= (n, n-1, \ldots, 1)$, so that
$\lambda=\mu+\delta$ is a strict partition of length $\ell(\lambda)=n$. Then
for $\a=(a_1,a_2,\ldots)$ we have:   
\begin{align}
P_\lambda(\x;\y|\a) &=\ds \prod_{1\leq i \leq n} x_i \ \prod_{1\leq i<j \leq n} (x_i + y_j)\ s_{\mu} (\x|\a)\,; \hfill\label{eqn-tok-fact-Psfn}\\ 
Q_\lambda(\x;\y|\a) &=\ds \prod_{1\leq i\leq j \leq n} (x_i + y_j)\ s_{\mu} (\x|\a)\,; \hfill\label{eqn-tok-fact-Qsfn} 
\end{align}
or, equivalently,
\begin{align}
   \ds  \sum_{P\in{\cal P}^\lambda(\n,\n')} \wgt(P)  &= \ds \prod_{1\leq i \leq n} x_i \ \prod_{1\leq i<j \leq n} (x_i + y_j) 
		              \  \sum_{T\in{\cal T}^\mu(\n)} \wgt(T)    \,; \hfill\label{eqn-tok-fact-PT} \\  
	\ds	 \sum_{Q\in{\cal Q}^\lambda(\n,\n')} \wgt(Q)  &= \ds  \prod_{1\leq i\leq j \leq n} (x_i + y_j) 
		              \  \sum_{T\in{\cal T}^\mu(\n)} \wgt(T)    \,;  \hfill\label{eqn-tok-fact-QT}
\end{align}
where 
\begin{equation}\label{eqn-Pwgt-Twgt}
   \wgt(P)=\!\!\!\prod_{(i,j)\in SF^\lambda}\!\!\!\wgt(\P_{ij}); ~~\wgt(Q)=\!\!\!\prod_{(i,j)\in SF^\lambda}\!\!\!\wgt(\Q_{ij});
   ~~\wgt(T)=\!\!\!\prod_{(i,j)\in F^\mu}\!\!\!\wgt(\T_{ij}),
\end{equation}	
with $\wgt(\P_{ij})$, $\wgt(\Q_{ij})$ and $\wgt(\T_{ij})$ given by
\begin{equation}\label{eqn-pwgt-twgt}
\begin{array}{|l|l|l|l|}
\hline
\P_{ii}&\wgt(\P_{ii})&\P_{ij}~~(i<j)&\wgt(\P_{ij})~~(i<j)\cr
\hline
 k&x_k &k&x_k+a_{j-i}\cr
   &   &k'&y_k-a_{j-i}\cr
\hline
\Q_{ii}&\wgt(\Q_{ii})&\Q_{ij}~~(i<j)&\wgt(\Q_{ij})~~(i<j)\cr
\hline
k&x_k&k&x_k+a_{j-i}\cr
k'&y_k&k'&y_k-a_{j-i}\cr
\hline
\end{array}
\qquad\hbox{and}\qquad
\begin{array}{|l|l|}
\hline
\T_{ij}&\wgt(\T_{ij}) \cr
\hline
 k&x_k+a_{k+j-i}\cr 
\hline
\end{array}
\end{equation}
\end{Theorem}


In specifying the weights as above, advantage has been taken of the fact that both $s_\lambda(\x\,|\,\a)$ and $Q_\lambda(\x;\y\,|\,\a)$ are independent 
of $a_0$ in our original definitions (\ref{eqn-factSchur}) and (\ref{eqn-factQSchur}), while $a_0$ is set equal to $0$ in the definition
of $P_\lambda(\x;\y\,|\,\a)$ in (\ref{eqn-factPSchur}). It might also be noted that under the hypothesis of this Theorem that $\ell(\lambda)=n$, 
the diagonal entries of any $S\in{\cal S}^\lambda([\n])$ are 
necessarily $1,2,\ldots,n$. It follows that the contributions of diagonal entries to every summand of $P_\lambda(\x;\y|\a)$ and to every summand of $Q_\lambda(\x;\y|\a)$ 
yield the factors $\prod_{i=1}^n x_i$ and $\prod_{i=1}^n (x_i+y_i)$, respectively. Since these factors represent the only difference between 
the expressions on the right hand sides of (\ref{eqn-tok-fact-Psfn}) and of (\ref{eqn-tok-fact-Qsfn}), in order to prove Theorem~\ref{thm-main}
it suffices only to prove the required results for either just $P_\lambda(\x;\y|\a)$ or just $Q_\lambda(\x;\y|\a)$. We choose to concentrate on the 
case $P_\lambda(\x;\y|\a)$ and construct the proof of (\ref{eqn-tok-fact-PT}).

In order to 
do this we make use of non-intersecting lattice path interpretations of the two sums 
appearing in (\ref{eqn-tok-fact-PT}), allowing each of them to be expressed in determinantal form by means of 
two lemmas, Lemma~\ref{lem-det-Twgt} and Lemma~\ref{lem-det-Pwgt} below, whose proofs we defer to the next section. 
A third, highly technical lemma, Lemma~\ref{lem-technical}, is then required
that allows us to proceed by way of simple row operations on the determinant representing the left hand side
of (\ref{eqn-tok-fact-PT}) to the required factorisation on the right. 
In view of its technical nature the proof of Lemma~\ref{lem-technical} is also deferred to the next section.

We begin with the determinantal expression for $s_\mu(\x|\a)$.
For the subsequence $\tilde{\x}=(x_k,x_{k+1},\ldots,x_n)$ of $\x$ with $1\leq k\leq n$ and $\a=(a_0,a_1,a_2,\ldots)$ 
let $h_m(\tilde{\x}|\a)=s_{(m)}(\tilde{\x}|\a)$ for all positive integers $m$. 
Then it follows from (\ref{eqn-factSchur}) that
\begin{equation}\label{eqn-hm}
\begin{array}{c}
     h_m(\tilde{\x}|\a) =  h_m(x_k,x_{k+1},x_{k+2},\ldots,x_n|\a)\cr\cr
\ds ~~~~~~~~~~~	= \sum_{k\leq i_1\leq i_2\leq \cdots\leq i_m\leq n}
			(x_{i_1}+a_{i_1-k+1})(x_{i_2}+a_{i_2-k+2})\cdots (x_{i_m}+a_{i_m-k+m})\,.
\end{array}			
\end{equation}

In terms of these single row factorial Schur functions we have the following determinantal identity 
that is originally due to Chen, Li and Louck~\cite{CLL}:
\begin{Lemma}\label{lem-det-Twgt}
Let $\mu$ be a partition of length $\ell(\mu)\leq n$, then 
\begin{equation}
   s_\mu(\x|\a)= \sum_{T\in{\cal T}^\lambda(\n)} \wgt(T) 
	     =  \det_{1\leq k,\ell\leq n} \big(\, h_{\mu_\ell-\ell+k}(x_k,x_{k+1},x_{k+2},\ldots,x_n|\a)\,\big)  \,,
   \label{eqn-det-Twgt}
\end{equation}
where $h_m(\tilde{\x}|\a)=1$ if $m=0$, and $h_m=0$ if $m<0$. %
\end{Lemma}

To set up the relevant determinantal expression for $P_\lambda(\x;\y|\a)$ we require certain shifted restricted versions 
$q_m(\tilde{\x};\tilde{\y}|\a)$ of the factorial generalised Schur $Q$ functions. 
Here shifts are associated with the introduction of an operator $\tau$~\cite{Macdonald} whose action
on $\a=(a_0,a_1,a_2,\ldots)$ is such that $\tau\a=(a_1,a_2,a_3,\ldots)$, so that in acting on any function of $\a$ each $a_i$ 
is replaced by $a_{i+1}$. For any $p,q$ and $r$ such that $1\leq p<q\leq r\leq n$ let $\tilde{\x}=(x_p,x_{p+1},\ldots,x_r)$ 
and $\tilde{\y}=(y_{q},y_{q+1},\ldots,y_r)$ be subsequences of our original sequences $\x$ and $\y$, respectively, and then let
\begin{equation}\label{eqn-qm}
  q_m(x_p,x_{p+1},\ldots,x_{q-1},y_{q},x_{q},y_{q+1},x_{q+1},\ldots,y_r,x_r|\a) = Q_{(m)}(\tilde{\x};\tilde{\y}|\tau\a)\,,
\end{equation}
where in evaluating the right hand side the entries in the one-rowed primed tableaux $P$ of $Q_{(m)}$ are taken from the 
alphabet $p<(p+1)<\cdots<(q-1)<q'<q<(q+1)'<(q+1)<\cdots<r'<r$	with repetitions allowed for unprimed entries but not for primed entries, 
and with $k'$ allowed in the box $(1,1)$ on the main diagonal if and only if $q\leq k\leq r$. The shift due to $\tau$ is 
such that an unprimed entry $k$ in column $j$ is weighted $x_k+a_j$ and a primed entry $k'$ in column $j$ is weighted $y_k-a_j$. Thus 
\begin{equation}\label{eqn-qm-z}
\begin{array}{c}
  q_m(x_p,x_{p+1},\ldots,x_{q-1},y_{q},x_{q},y_{q+1},x_{q+1},\ldots,y_r,x_r|\a) \cr\cr
	\ds ~~~~~~~~~~~ = \sum_{p\leq i_1\leq i_2\leq\cdots\leq i_m\leq r} \sum_\z (z_{i_1}\pm a_1)(z_{i_2}\pm a_2)\cdots(z_{i_m}\pm a_m)\,,
\end{array}				
\end{equation}
where the sum over $\z$ allows factors $(z_k\pm a_j)=(x_k+a_j)$ or $(y_k-a_j)$ to appear according as $z_k=x_k$ or $y_k$,
with several factors of the form $(x_k+a_j)(x_k+a_{j+1})\cdots$ allowed for any $k$ with $p\leq k\leq r$
but at most one factor $(y_k-a_j)$ allowed for any $k$ such that $q\leq k\leq r$, and no others.

This allows us to express the left hand side of (\ref{eqn-tok-fact-PT}) in the form of a determinant by means
of the following key lemma: 
\begin{Lemma}\label{lem-det-Pwgt}
Let $\lambda$  
be a strict partition of length $\ell(\lambda)=n$.  
Then we have
\begin{equation}
   \sum_{P\in{\cal P}^\lambda(\n,\n')} \wgt(P) =  \det_{1\leq k,\ell\leq n} 
	         \big(\, x_k q_{\lambda_{\ell}-1}(x_k, y_{k+1}, x_{k+1}, y_{k+2}, \ldots, y_n, x_n|\a)\,\big)  
   \label{eqn-det-Pwgt}
\end{equation}
\end{Lemma}

The evaluation of this determinant may be accomplished by way of a technical lemma. In order to state this it is necessary 
to introduce a second type of shift operator $\sh$ that unlike $\tau$ is linked to
letters of the alphabet. The action of $\sh$ inserted in the $j$th position in 
$q_m(z_1,z_2,\ldots,z_r|\a)$ gives $q_m(z_1,z_2,\ldots,z_{j-1},\sh z_j,z_{j+1},\ldots,z_r|\a)$ in which
every linear factor $(z_i+a_s)$ or $(z_i-a_t)$ of $q_m(z_1,z_2,\ldots,z_r|\a)$ is replaced by 
$(z_i+a_{s+1})$ or $(z_i-a_{t+1})$, respectively, if and only if $i>j$.
In other words the insertion of the operator $\sh$ increases by $1$ the index of $a$ in every linear factor associated with each parameter to its right.
Repeated insertions of shift operators $\sh$ are allowed, either at the same or at different points. Powers such as $\sh^p$ inserted at a single point
increase by $p$ the index of $a$ in every linear factor associated with each parameter to its right. 
Now we can state our technical lemma.

\begin{Lemma}\label{lem-technical}
For all $m\geq 1$ we have:  \\
\indent  a) For $1\leq p<n$ 
\begin{equation}\label{eqn-lem-a}
\begin{array}{l}
q_m(x_p, y_{p+1}, x_{p+1}, y_{p+2}, x_{p+2}, \ldots, y_n, x_n|\a)
 - q_m(x_{p+1}, y_{p+2}, x_{p+2}, \ldots, y_n, x_n|\a)\\ 
= (x_p+y_{p+1}) q_{m-1}(x_p,\sh x_{p+1}, y_{p+2}, x_{p+2}, \ldots, y_n, x_n|\a)
\end{array}
\end{equation}

and more generally,

b) For $1\leq p<q\leq n$ 
\begin{equation}\label{eqn-lem-b}
\begin{array}{c}
q_m(x_p, \sh x_{p+1}, \sh x_{p+2}, \cdots, \sh x_{q-1}, y_q, x_q, y_{q+1}, x_{q+1}, \ldots, y_n, x_n|\a) \\
 - q_m(x_{p+1}, \sh x_{p+2}, \cdots, \sh x_{q-1}, \sh  x_q, y_{q+1}, x_{q+1}, \ldots, y_n, x_n|\a)\\ 
= (x_p+y_q) q_{m-1} (x_p, \sh x_{p+1}, \cdots, \sh x_{q-1}, \sh x_q, y_{q+1}, x_{q+1}, \ldots, y_n, x_n|\a)
\end{array}
\end{equation}
\end{Lemma}

Given these three lemmas we have enough to prove our main Theorem \ref{thm-main}.

{\bf Proof of Theorem \ref{thm-main}:}  Lemma \ref{lem-det-Pwgt} expresses the left hand side of (\ref{eqn-tok-fact-PT}) 
as a determinant from which we can extract $x_k$ from each row for $k=1,2,\ldots,n$ to give 
\begin{equation}
   \sum_{P\in{\cal P}^\lambda(\n,\n')} \wgt(P) =  \prod_{i=1}^n x_i\ \ \det_{1\leq k,\ell\leq n} 
	         \left(q_{\lambda_{\ell}-1}(x_k, y_{k+1}, x_{k+1}, y_{k+2}, \ldots, y_n, x_n|\a)\,\right)  \,.
\end{equation}
Subtracting row $k+1$ from row $k$ of the determinant for $k=1,2,\ldots,n-1$ gives a new determinant in which 
the $(k,\ell)$th elements is given by left hand side of (\ref{eqn-lem-a}) with $p=k$, $q=k+1$ and $m=\lambda_{\ell-1}$, while the $n$th row
remains unaltered with elements $q_{\lambda_{\ell-1}}(x_n|\a)$. 
Applying (\ref{eqn-lem-a}) and extracting a common factor of $(x_k+y_{k+1})$ from the $k$th row then gives
\begin{equation}
\begin{array}{l}
\ds  \sum_{P\in{\cal P}^\lambda(\n,\n')} \wgt(P) = \prod_{i=1}^n x_i\ \prod_{i=1}^{n-1} (x_i+y_{i+1})\cr
\ds 	~~~~~~~~~~~~~ \times\  \det_{1\leq k,\ell\leq n} \  
	         \left( \begin{array}{c}
					q_{\lambda_{\ell}-2}(x_k, \sh x_{k+1}, y_{k+2}, x_{k+2} \ldots, y_n, x_n|\a)\cr
					q_{\lambda_{\ell}-1}(x_n|\a)
					       \end{array}
					\,\right)  \,,
\end{array}					
\end{equation}
where we have distinguished between elements in the first $n-1$ rows and the last row.
 
We can then use the same procedure of subtracting row $k+1$ from row $k$ of the above determinant,
this time for $k=1,2,\ldots,n-2$ to give a new determinant in which the $(k,\ell)$th element is given by the 
left hand side of part b) of Lemma~\ref{lem-technical} with $p=k$, $q=k+2$ and $m=\lambda_{\ell-2}$. Applying
(\ref{eqn-lem-b}) and extracting a common factor of $(x_k+y_{k+2})$ from the $k$th row then gives
\begin{equation}
\begin{array}{l}
\ds   \sum_{P\in{\cal P}^\lambda(\n,\n')} \wgt(P) = \prod_{i=1}^n x_i\  \prod_{i=1}^{n-1} (x_i+y_{i+1})\  \prod_{i}^{n-2} (x_i+y_{i+2})\ \cr 
\ds ~~~~~~~~~~~~~ \times \ \det_{1\leq k,\ell\leq n} \left( \begin{array}{c}
					q_{\lambda_{\ell}-3}(x_k, \sh x_{k+1}, \sh x_{k+2}, y_{k+3} , x_{k+3},\ldots, y_n, x_n|\a)\cr
					q_{\lambda_{\ell}-2}(x_{n-1}, \sh x_n|\a)\cr
					q_{\lambda_{\ell}-1}(x_n|\a)
					       \end{array}
					\,\right)  \,,
\end{array}					
\end{equation}
where this time we have distinguished between elements in the first $n-2$ rows and the last $2$ rows.

Continuing in this way we obtain
\begin{eqnarray}
\lefteqn{\sum_{P\in{\cal P}^\lambda(\n,\n')}\!\!\! \wgt(P)}\\
& =& \prod_{i=1}^n x_i\  \prod_{1\leq i<j\leq n}\!\!(x_i+y_j)
\ \det_{1\leq k,\ell\leq n} \left( 
					q_{\lambda_{\ell}-n+k-1}(x_k, \sh x_{k+1}, \sh x_{k+2},\ldots, \sh x_n|\a)
									\,\right).
\end{eqnarray}
However  
\begin{equation}
\begin{array}{l}
  q_m(x_k, \sh x_{k+1}, \sh x_{k+2},\ldots, \sh x_n|\a) \cr
\ds ~~~~~~~~	=\sum_{k\leq i_1\leq i_2\leq \cdots\leq i_m \leq n} (x_{i_1}+a_{i_1-k+1})(x_{i_2}+a_{i_2-k+2})\cdots(x_{i_m}+a_{i_m-k+m})\,,
\end{array}
\end{equation}
where account has been taken of the fact that there are precisely $(i_j-k)$ shift operators
$\sh$ to the left of $x_{i_j}$ in the argument of $q_m$. This will be recognised as coinciding with the definition
of $h_m(x_k,x_{k+1},\ldots,x_n|\a)$ given in (\ref{eqn-hm}).
Then the use of Lemma~\ref{lem-det-Twgt} with $\mu_\ell=\lambda_\ell-n+\ell-1$ completes the proof of 
(\ref{eqn-tok-fact-PT}) and thereby that of Theorem~\ref{thm-main}.
It remains only to prove the validity of our three lemmas, Lemmas~\ref{lem-det-Twgt},~\ref{lem-det-Pwgt} and \ref{lem-technical}.
\qed
\bigskip

\section{Proofs of Lemmas~\ref{lem-det-Twgt},~\ref{lem-det-Pwgt} and \ref{lem-technical}}
\label{sec:lem-proofs}

In each case we follow the lattice path approach of Okada \cite{Okadapartial}, employing a variation on the usual 
Gessel-Viennot-Lindstr\"om argument (see in particular Okada \cite{Okadapartial} and Stembridge \cite{Stembridge}).  
In the case of Lemma~\ref{lem-det-Twgt} a similar proof by way of a lattice path interpretation
has been offered by Chen, Li and Louck~\cite{CLL}, but we offer an independent lattice path proof here that takes 
particular advantage of the precise form of $h_m(\tilde{\x}\,|\,\a)$ given in (\ref{eqn-hm}), since it is this form
that we have just seen emerging in a natural way in the application of Lemma~\ref{lem-technical} to the proof of
Theorem~\ref{thm-main}. Moreover, it is our lattice path proof of Lemma~\ref{lem-det-Twgt} that sets the scene for 
our rather similar lattice path proof of Lemma~\ref{lem-det-Pwgt}.
\bigskip

\noindent{\bf Proof of Lemma~\ref{lem-det-Twgt}:}
We adopt matrix coordinates $(i,j)$ for  
lattice points with $i=1,2,\ldots,n$ specifying row labels from
top to bottom, and $j=1,2,\ldots,\mu_1+n$ specifying column labels from left to right. Each lattice path that we are interested in 
is a continuous path from some 
$P_i=(i,n-i+1)$ with $i\in\{1,2,\ldots,n\}$ to some $Q_j=(n+1,j)$
with $j\in\{\mu_1+n,\mu_2+n-1,\ldots,\mu_n+1\}$. Such a path consists of a sequence of horizontal or vertical edges with the last edge vertical. 

\begin{figure*}[htbp]
\begin{center}
\begin{equation*}
T\ = \ 
\YT{0.2in}{}{
 {1,2,4},
 {2,3},
 {4,4},
 {5}
}
\qquad 
\vcenter{\hbox{
\begin{tikzpicture}[x={(0in,-0.25in)},y={(0.25in,0in)}] 
\foreach \j in {5,...,8} \draw(1,\j)node{$\bullet$};
\foreach \j in {4,...,8} \draw(2,\j)node{$\bullet$};
\foreach \j in {3,...,8} \draw(3,\j)node{$\bullet$};
\foreach \j in {2,...,8} \draw(4,\j)node{$\bullet$};
\foreach \j in {1,...,8} \draw(5,\j)node{$\bullet$};
\draw[->](0+0.2,7-0.2)to(1-0.2,6+0.2);
\draw[->](0+0.2,8-0.2)to(1-0.2,7+0.2);
\draw[->](0+0.2,9-0.2)to(1-0.2,8+0.2);
\draw[->](1+0.2,9-0.2)to(2-0.2,8+0.2);
\draw[->](1+0.2,9-0.2)to(2-0.2,8+0.2);
\draw[->](2+0.2,9-0.2)to(3-0.2,8+0.2);
\draw[->](3+0.2,9-0.2)to(4-0.2,8+0.2);
\draw[->](4+0.2,9-0.2)to(5-0.2,8+0.2);
\draw(1,4.5)node{$P_1$};
\draw(2,3.5)node{$P_2$};
\draw(3,2.5)node{$P_3$};
\draw(4,1.5)node{$P_4$};
\draw(5,0.5)node{$P_5$};
\draw(7,1)node{$Q_1$};
\draw(7,3)node{$Q_3$};
\draw(7,5)node{$Q_5$};
\draw(7,6)node{$Q_6$};
\draw(7,8)node{$Q_8$};
\draw(0,7)node{$a_1$};
\draw(0,8)node{$a_2$};
\draw(0,9)node{$a_3$};
\draw(1,9)node{$a_4$};
\draw(2,9)node{$a_5$};
\draw(3,9)node{$a_6$};
\draw(4,9)node{$a_7$};
\draw[draw=black,thick] (1,5)to(1,6);
\draw[draw=black,thick] (1,6)to(2,6);
\draw[draw=black,thick] (2,6)to(2,7);
\draw[draw=black,thick] (2,7)to(4,7);
\draw[draw=black,thick] (4,7)to(4,8);
\draw[draw=black,thick] (4,8)to(6,8);
\draw[draw=black,thick] (2,4)to(2,5);
\draw[draw=black,thick] (2,5)to(3,5);
\draw[draw=black,thick] (3,5)to(3,6);
\draw[draw=black,thick] (3,6)to(6,6);
\draw[draw=black,thick] (3,3)to(4,3);
\draw[draw=black,thick] (4,3)to(4,5);
\draw[draw=black,thick] (4,5)to(6,5);
\draw[draw=black,thick] (4,2)to(5,2);
\draw[draw=black,thick] (5,2)to(5,3);
\draw[draw=black,thick] (5,3)to(6,3);
\draw[draw=black,thick] (5,1)to(6,1);
\end{tikzpicture}
}}
\end{equation*}
\end{center}
\caption{Example of  the lattice paths for a given semistandard tableau.}
\label{fig-TtoLP}
\end{figure*}
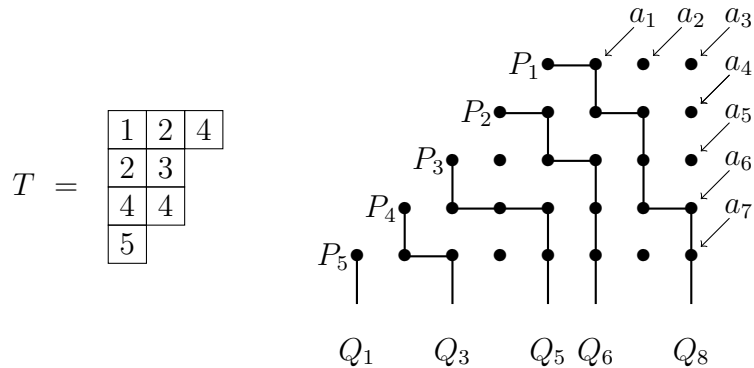

Each semistandard tableau $T$ of shape $\mu$ defines a set of non-intersecting lattice paths, one for each row of $T$.
The path associated with the $i$th row of $T$ starts at $P_i$ and ends at $Q_j$ with $j=\mu_i+n-i+1$. 
On this path each entry $k$ in the $\ell$th column of $T$ gives rise to a horizontal edge from $(k,j-1)$ to $(k,j)$ 
with $j=n-k+\ell$, and vertical edges are added so as to make the path continuous. It is easy to see from the properties T1-T3 of 
Section~\ref{sec:background} that the paths are non-intersecting. This is exemplified in Figure~\ref{fig-TtoLP} in the case $\mu=(3,2,2,1,0)$ 
and $T$ as given in (\ref{eqn-TSP}).

The map we have described from $T$ to a set $L$ of non-intersecting lattice paths is a bijection as can be seen by reversing the argument and mapping consecutive horizontal edges at level $k$ along a path starting at $P_i$ to entries $k$ in the $i$th 
row of $T$. The non-intersecting nature of the paths ensures that $T$ constructed in this way is a semistandard
tableau as required.
  
In order to recover $\wgt(T)$ as defined through (\ref{eqn-pwgt-twgt}) from the set of lattice paths it is important to 
note that the entries $k$ of $T$ are associated with the $k$th row of the lattice, and that the $\ell$th column of $T$ 
is associated with the $\ell$th diagonal of the lattice along which $k+j=n+1+\ell$. Each horizontal edge from $(k,j-1)$ to $(k,j)$ is weighted $x_k+a_{k+j-n-1}$ and each vertical edge is weighted $1$. With these asignments it follows that $\wgt(T)$ is just the product over all edges of these edge weights. 
Thus the left hand side of (\ref{eqn-det-Twgt}) is evaluated by summing over all sets of non-intersecting paths with the given 
end points $P_i$ and $Q_{\mu_i+n-i+1}$ 
with $i=1,2,\ldots,n$.

More generally, the total weight of all possible continuous lattice paths from $P_i$ to $Q_j$ by means of horizontal and vertical edges is given by some summand of $h_{m}(x_k,x_{k+1},\ldots,x_n|\a)$ with $m=i+j-n-1$. 
Then the usual argument~\cite{Okadapartial}, extended so as to allow a fixed set of end points determined as in our case by 
$\mu$, shows that the total weight of the set of all (intersecting and non-intersecting) lattice paths from the given set of starting points $P_i$ to the ending points $Q_j$, summed over all permutations of $Q_j$, is exactly the determinant of the matrix whose $(k,\ell)$th entry is $h_{\mu_\ell}(x_k, x_{k+1}, \ldots, x_n|\a)$, as required to complete the proof of 
Lemma~\ref{lem-det-Twgt}.  \qed
\bigskip

\noindent{\bf Proof of Lemma~\ref{lem-det-Pwgt}}

It is again convenient to adopt matrix coordinates $(i,j)$ for the lattice points with $i=1,2,\ldots,n$ specifying row labels from
top to bottom, and $j=1,2,\ldots,\lambda_1$ specifying column labels from left to right. This time each lattice path that we are interested 
in is a continuous path from some $P_i=(i,0)$ 
with $i\in\{1,2,\ldots,n\}$ to some $Q_j=(n+1,j)$ with 
$j\in\{\lambda_1,\lambda_2,\ldots,\lambda_n\}$. Such a path now consists of a sequence of horizontal, diagonal or vertical edges with the first edge horizontal and the last edge vertical. 

\begin{figure*}[htbp]
\begin{center}
\begin{equation*}
P\ = \
\SYT{0.2in}{}{
 {1,1,2',2,3',4},
   {2,3',3,3},
     {3,4',4},
       {4}
}
\qquad 
\vcenter{\hbox{
\begin{tikzpicture}[x={(0in,-0.25in)},y={(0.25in,0in)}] 
\foreach \i in {1,...,4} \foreach \j in {1,...,6} \draw(\i,\j)node{$\bullet$};
\draw(1,-1)node{$P_1$};
\draw(2,-1)node{$P_2$};
\draw(3,-1)node{$P_3$};
\draw(4,-1)node{$P_4$};
\draw(6,1)node{$Q_1$};
\draw(6,3)node{$Q_3$};
\draw(6,4)node{$Q_4$};
\draw(6,6)node{$Q_6$};
\foreach \j in {3,...,7} \draw[->](0+0.2,\j-1)to(1-0.3,\j-1);
\foreach \j in {1,...,5} \draw(0,\j+1)node{$a_{\j}$};
\foreach \i in {1,...,4} \draw[draw=black,thick] (\i,0)to(\i,1);
\draw[draw=black,thick] (1,1)to(1,2);
\draw[draw=black,thick] (1,2)to(2,3);
\draw[draw=black,thick] (2,3)to(2,4);
\draw[draw=black,thick] (2,4)to(3,5);
\draw[draw=black,thick] (3,5)to(4,5);
\draw[draw=black,thick] (4,5)to(4,6);
\draw[draw=black,thick] (4,6)to(5,6);
\draw[draw=black,thick] (2,1)to(3,2);
\draw[draw=black,thick] (3,2)to(3,4);
\draw[draw=black,thick] (3,4)to(5,4);
\draw[draw=black,thick] (3,1)to(4,2);
\draw[draw=black,thick] (4,2)to(4,3);
\draw[draw=black,thick] (4,3)to(5,3);
\draw[draw=black,thick] (4,1)to(5,1);
\end{tikzpicture}
}}
\end{equation*}
\end{center}
\caption{Example of  the lattice paths for a given primed shifted tableau.}
\label{fig-PtoLP}
\end{figure*}
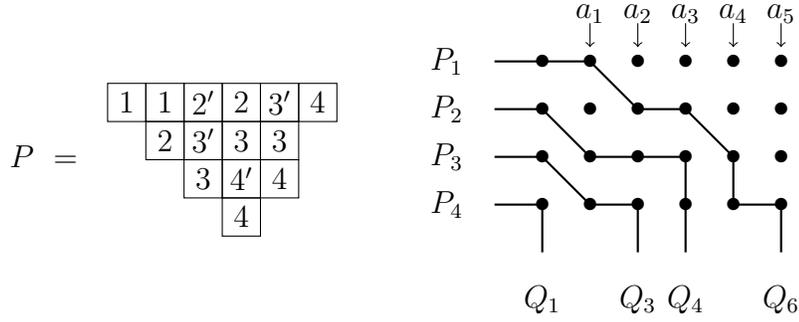

Each primed shifted tableau $P$ of shape $\lambda$ defines a set of non-intersecting lattice paths, one for each row of $P$.
The path associated with the $i$th row of $P$ starts at $P_i$ and ends at $Q_j$ with $j=\lambda_i$. 
Each unprimed entry $k$ in the $j$th diagonal of $P$ gives rise to a horizontal edge from $(k,j-1)$ to $(k,j)$ and each primed entry $k'$ in the $j$th diagonal of $P$ gives rise to a diagonal edge from $(k-1,j-1)$ to $(k,j)$ with vertical edges being added so as to make the path continuous. It is easy to see from the properties P1-P5 of Section~\ref{sec:background} that the paths are non-intersecting. This is exemplified in Figure~\ref{fig-PtoLP} in the case $\lambda=(6,4,3,1)$ and $P$ as given in (\ref{eqn-TSP}). 

The map we have described from $P$ to a set $L$ of non-intersecting lattice paths is a bijection as can be seen by reversing the argument and mapping consecutive horizontal and diagonal edges along a path starting at $P_i$ to entries $k$ and $k'$, respectively, in the $i$th row of $P$. The non-intersecting nature of the paths ensures that $P$ constructed in this way is a 
primed shifted tableau as required.
  
In order to recover $\wgt(P)$ as defined through (\ref{eqn-pwgt-twgt}) from the set of lattice paths it is important to 
note that the entries $k$ and $k'$ in $P$ are associated with the $k$th row of the lattice, and that the $\ell$th diagonal of $P$ is associated with the $\ell$th column of the lattice. 
Since $\ell(\lambda)=n$ and the $i$th row of $P$ necessarily begins with an unprimed entry $i$, the first horizontal edge of the path starting at $P_i$ is weighted $x_i$. As far as the remaining edges of the set of lattice paths is concerned, any
horizontal edge from $(k,\ell)$ to $(k,\ell)$ is weighted $x_k+a_{\ell-1}$, any diagonal edge from $(k-1,\ell-1)$ to $(k,\ell)$
is weighted $y_k-a_{\ell-1}$ and each vertical edge is weighted $1$. With these asignments it follows that $\wgt(P)$ is just the product over all edges of these edge weights. Thus the left hand side of (\ref{eqn-det-Pwgt}) is evaluated by summing over all sets of non-intersecting paths with the given starting points $P_i$ and $Q_{\lambda_i}$ with $i=1,2,\ldots,n$.

Given this framework, it is not hard to see that the total weight of all continuous lattice paths from $P_k$ to $Q_\ell$ by means of the three types of edge, horizontal, diagonal and vertical is given by  
$x_k q_{\lambda_{\ell}-1}(x_k, y_{k+1}, x_{k+1}, \ldots, y_n, x_n|\a)$. 
Then Okada's argument in \cite{Okadapartial}, extended so as to allow a fixed set of end points determined as in our case by 
$\lambda$, shows that the total weight of the set of all (intersecting and non-intersecting) lattice paths from the given set of starting points $P_i$ to the ending points $Q_j$, summed over all permutations of $Q_j$, is exactly the determinant of the matrix whose $(k,\ell)$th entry is $x_k q_{\lambda_{\ell}-1}(x_k, y_{k+1}, x_{k+1}, \ldots, y_n, x_n|\a)$, as required to complete the proof of Lemma~\ref{lem-det-Pwgt}.  \qed
\bigskip

Now all that remains is to prove our technical Lemma~\ref{lem-technical}. Although part a) is really a special 
case of part b) its proof is somewhat different from the special case, so we treat it seperately as follows.

\noindent{\bf Proof of Lemma~\ref{lem-technical} part a):}~~The identity to be proved takes the form
\begin{equation}\label{eqn-lem-a-proof}
\begin{array}{l}
q_m(x_p, y_{p+1}, x_{p+1}, y_{p+2}, x_{p+2}, \ldots, y_n, x_n|\a)
 - q_m(x_{p+1}, y_{p+2}, x_{p+2}, \ldots, y_n, x_n|\a)\\ 
= (x_p+y_{p+1}) q_{m-1}(x_p,\sh x_{p+1}, y_{p+2}, x_{p+2}, \ldots, y_n, x_n|\a)\,.
\end{array}
\end{equation}
Any term in $q_m(x_p, y_{p+1}, x_{p+1}, y_{p+2}, x_{p+2}, \ldots, y_n, x_n|\a)$ that does 
not contain either $x_p$ or $y_{p+1}$ must be of $m$th degree in variables chosen from the set 
$x_{p+1}, y_{p+2}, x_{p+2}, \ldots, y_n,$ $x_n$, that is to say a term in 
$q_m(x_{p+1}, y_{p+2}, x_{p+2}, \ldots, y_n, x_n|\a)$.  
Thus these terms cancel between the expressions in the first and second lines of (\ref{eqn-lem-a-proof}).

Now the remaining terms on the left either contain no $y_{p+1}$ but at least one $x_{p}$ and therefore
constitute partial sums of the form 
\[
q_{k-1}(x_p|\a)\ (x_p+a_k)\ q_{m-k}(\sh^{k}x_{p+1}, y_{p+2}, x_{p+2}, \ldots, y_n, x_n|\a)
\]
for some $k$, or they contain one, but no more than one $y_{p+1}$, together with an arbitrary number of $x_p$'s, and therefore 
constitute partial sums of the form
\[
q_{k-1}(x_p|\a)\ (y_{p+1}-a_k)\ q_{m-k}(\sh^{k}x_{p+1}, y_{p+2}, x_{p+2}, \ldots, y_n, x_n|\a)
\]
for some $k$. Pairing them in this way and adding gives
\[
\begin{array}{c}
\ds \sum_{k=1}^m (x_p+a_k +y_{p+1} - a_k)\ q_{k-1} (x_p|\a)\ q_{m-k}(\sh^{k}x_{p+1}, y_{p+2}, x_{p+2}, \ldots, y_n, x_n|\a) \\
~~ = ~~ (x_p+y_{p+1} )\ q_{m-1} (x_p, \sh x_{p+1}, y_{p+2}, x_{p+2}, \ldots, y_n, x_n|\a) \,.
\end{array}
\]
\qed
\bigskip

\noindent{\bf Proof of Lemma~\ref{lem-technical} part b):}~~This time the identity to be proved takes the form:
\begin{equation}\label{eqn-lem-b-proof}
\begin{array}{c}
q_m(x_p, \sh x_{p+1}, \sh x_{p+2}, \cdots, \sh x_{q-1}, y_q, x_q, y_{q+1}, x_{q+1}, \ldots, y_n, x_n|\a) \\
 - q_m(x_{p+1}, \sh x_{p+2}, \cdots, \sh x_{q-1}, \sh  x_q, y_{q+1}, x_{q+1}, \ldots, y_n, x_n|\a)\\ 
= (x_p+y_q) q_{m-1} (x_p, \sh x_{p+1}, \cdots, \sh x_{q-1}, \sh x_q, y_{q+1}, x_{q+1}, \ldots, y_n, x_n|\a)\,.
\end{array}
\end{equation}
We proceed by showing that we can cluster the terms arising from the first line of (\ref{eqn-lem-b-proof}) in a particular way,
then reorganise each cluster so as to separate off terms containing the required factor $(x_p+y_q)$ contributing to the right hand side 
appearing on the third line of (\ref{eqn-lem-b-proof}),
leaving terms we call children that contribute to subsequent clusters, with a final remainder that can be identified with the second line of (\ref{eqn-lem-b-proof}).
Although, strictly speaking, what we construct is not a tree, there is a hierarchical structure that is similar to that of a tree, and we use some tree terminology in what follows. Each cluster is associated with a summand $u_{k-1}$ appearing in the expansion of $q_{k-1}(x_p, \sh x_{p+1}, \sh x_{p+2}, \cdots, \sh x_{q-1}|\a)$.
These clusters for fixed $k=1,2,\ldots,m$ are arranged in reverse lexicographic order from top to bottom to form the trunk of an (upside down) tree.

We begin by identifying the topmost ``root cluster'' of the tree corresponding to $k=m$ 
as the sum of terms differing in content in only one location, namely the last, that is given by:
\begin{eqnarray*}
(x_p + a_{1})(x_p +a_{2}) \ldots(x_p+a_{m-1}) (x_p+a_{m})+\\
(x_p + a_{1})(x_p +a_{2}) \ldots(x_p+a_{m-1}) (x_{p+1}+a_{m+1})+\\
(x_p + a_{1})(x_p +a_{2}) \ldots(x_p+a_{m-1}) (x_{p+2}+a_{m+2})+\\
\ldots~~~~+\\
(x_p + a_{1})(x_p +a_{2}) \ldots(x_p+a_{m-1}) (x_{q-1}+a_{m+q-p-1})+\\
(x_p + a_{1})(x_p +a_{2}) \ldots(x_p+a_{m-1}) (y_q-a_{m+q-p-1}) \,.
\end{eqnarray*}
This can be rewritten as
\begin{equation}
\resizebox{.99 \textwidth}{!}
{  
$
\begin{array}{l}
(x_p + a_{1})(x_p +a_{2}) \ldots(x_p+a_{m-1}) \times \cr
(x_p+a_{m} +x_{p+1}+a_{m+1} +x_{p+2}+a_{m+2} + 
\ldots + x_{q-1}+a_{m+q-p-1} + y_q-a_{m+q-p-1})\,.
 \end{array}
 $
  }
\end{equation}
By cyclically shifting the $a$ terms by one step and cancelling the terms $\pm a_{m+q-p-1}$, this can further be rewritten as
\begin{eqnarray*} 
(x_p + a_{1})(x_p +a_{2}) \ldots(x_p+a_{m-1}) (x_p)+\\  \nonumber
(x_p + a_{1})(x_p +a_{2}) \ldots(x_p+a_{m-1}) (x_{p+1}+a_{m})+\\  \nonumber
(x_p + a_{1})(x_p +a_{2}) \ldots(x_p+a_{m-1}) (x_{p+2}+a_{m+1})+\\  \nonumber
\ldots~~~~+\\ \nonumber
(x_p + a_{1})(x_p +a_{2}) \ldots(x_p+a_{m-1}) (x_{q-1}+a_{m+q-p-2})+\\  
(x_p + a_{1})(x_p +a_{2}) \ldots(x_p+a_{m-1}) (y_q) \,.
\end{eqnarray*}
The term
\[
(x_p+y_q)(x_p + a_{1})(x_p +a_{2}) \ldots(x_p+a_{m-1}) 
\]
can then be removed and is exactly the type of term required on the right hand side of (\ref{eqn-lem-b-proof}).  The remaining terms are what we call ``children''.  Each child is grouped with $q-p$ other terms (some of which are children of other terms, and some of which are original terms from the first line of (\ref{eqn-lem-b-proof})) such that the $q-p+1$ terms 
form a new cluster having the same form as the root cluster, namely the terms are identical in all but one location and in that one location they run through all possibilities from $x_p, x_{p+1}, \ldots, x_{q-1}, y_q$.  Given this, we can perform the same operations as in the case of the root cluster. That is to say adding these terms together, cyclically shifting the $a$ contributions, cancelling a pair of terms of the form $\pm a$ and separating out the terms containing the common factor $(x_p + y_q)$.

Iterating this procedure, we are left with terms we call ``leaves'' that contain neither $x_p$ nor $y_q$ and that 
exactly match the negative terms in the second line of (\ref{eqn-lem-b-proof}) and thus cancel out.

In order to be more precise it is necessary to define what we mean by a ``cluster''. For fixed $1\leq p<q\leq n$ and fixed $m$
each cluster is a sum of $q-p+1$ terms taking the form $u_{k-1}\ v_1\ w_{m-k}$ where
\begin{equation}\label{eqn-cluster}
\begin{array}{ccll}
u_{k-1}&=&(x_{i_1}+a_{j_1})(x_{i_2}+a_{j_2})\cdots(x_{i_{k-1}}+a_{j_{k-1}})\cr 
v_{1}&=&(x_p+a_{c_{p}}) + (x_{p+1}+a_{c_{p+1}}) + \cdots + (x_{q-1}+a_{c_{q-1}}) + (y_q-a_{c_{q-1}})     \cr
w_{m-k}&=&(z_{i_{k+1}}\pm a_{j_{k+1}})(z_{i_{k+2}}\pm a_{j_{k+2}})\cdots(z_{i_{m}}\pm a_{j_{m}})
\end{array}
\end{equation}
with $1\leq p\leq i_1\leq i_2\leq \cdots\leq i_{k-1}<q\leq i_{k+1}\leq i_{k+2}\leq \cdots\leq i_m\leq n$. 
Here $(z_i\pm a_j)=(x_i+a_j)$ or $(y_j-a_j)$ 
according as $z_i=x_i$ or $y_i$, respectively, with at most one term $(y_i-a_j)$ allowed for each $i>q$.
The indices on the factorial shifts $a$ are specified as follows:
\begin{equation}
\begin{array}{rcll}
     (x_{i_\ell}\+a_{j_\ell}) &=& (x_r\+a_{r-p+\ell})  &\hbox{if $i_\ell=r$ for $\ell=1,2,\ldots,k-1$}\cr
		 (x_{r}\+a_{c_{r}}) &=& (x_r\+a_{r-p+\ell}) &\hbox{where $(x_t\+a_{t-p+\ell})$ is the leftmost factor }\cr
		 & & & \hbox{of $u_{k-1}$with $t>r$} \cr
		 (z_{i_\ell}\!\pm\!a_{j_\ell}) &=& (z_s\!\pm\!a_{q-p+\ell-1})  &\hbox{if $i_\ell=s$ for $\ell=k+1,k+2,\ldots,m$}\cr
\end{array}    
\end{equation}

These conditions ensure that $u_{k-1}$ is a term in $q_{k-1}(x_p, \sh x_{p+1}, \sh x_{p+2}, \cdots, \sh x_{q-1}|\a)$
and that $w_{m-k}$ is a term in $q_{m-k}(\sh^{q-p+1} x_q, y_{q+1}, x_{q+1}, \ldots, y_n, x_n|\a)$,
while those of the form $u_{k-1}(x_r+a_{c_{r}})w_{m-k}$ with $r\geq i_{k-1}$ are standard in the sense that they appear in 
$q_m(x_p, \sh x_{p+1}, \sh x_{p+2}, \cdots, \sh x_{q-1}, y_q, x_q, y_{q+1}, x_{q+1}, \ldots, y_n, x_n|\a)$,
as is also true of $u_{k-1}(y_q-a_{c_{q-1}})w_{m-k}$. The remaining terms of the form $u_{k-1}(x_r+a_{c_{r}})w_{m-k}$ with $r<i_{k-1}$
are non-standard. It is these non-standard terms that are children of terms in ancestral clusters.
For $r<i_{k-1}$ the indices $c_{r}$ have been chosen so that if the leftmost factor $(x_t+a_{s-p+\ell})$ of $u_{k-1}$ with $t>r$ is replaced 
by $(x_r+a_{r-p+\ell})$ then the result $u'_{k-1}$ is still a term in $q_{k-1}(x_p, \sh x_{p+1}, \sh x_{p+2}, \cdots, \sh x_{q-1}|\a)$
with $u'_{k-1}$ preceding $u_{k-1}$ in reverse lexicographic order, and $u'_{k-1}(x_t+a_{t-p+\ell})w_{m-k}$, when reverse cyclically shifted, constituting
a legitimate term in what we call an ancestral cluster.

This is exemplified in the case $p=3$, $q=9$, $n=10$, $k=7$ and $m=11$ by the cluster:
\begin{equation}\label{eqn-cluster-ex}
\begin{array}{c}
(x_3\!+\!a_1)(x_3\!+\!a_2)(x_4\!+\!a_4)(x_6\!+\!a_7)(x_6\!+\!a_8)(x_6\!+\!a_9)\cr
\times((\red{x_3\!+\!a_3})+(\red{x_4\!+\!a_5})+(\red{x_5\!+\!a_6})+(\blue{x_6\!+\!a_{10}})+(\blue{x_7\!+\!a_{11}})+(\blue{x_8\!+\!a_{12}})+(\blue{y_{9}\!-\!a_{12}}))\cr
\times(x_9\!+\!a_{13})(x_9\!+\!a_{14})(y_{10}\!-\!a_{15})(x_{10}\!+\!a_{16})
\end{array}
\end{equation}
where the terms in blue are standard and those in red are non-standard. If we adopt the notation
$s_t=(x_s+a_t)$ and $s'_t=(y_s-a_t)$ this cyclically shifted cluster may be represented by
\begin{equation}\label{eqn-cluster-ex-short}
3_1 3_2 4_4 6_7 6_8 6_9 (\red{3_3\+4_5\+5_6}\+\blue{6_{10}\+7_{11}\+8_{12}\+9'_{12}}) 9_{13} 9_{14} 10'_{15} 10_{16}\,.
\end{equation}

Returning to the general case, cyclically shifting the $a$'s and cancelling the pair $\pm a$
in (\ref{eqn-cluster}) we can rewrite $v_1$ as
\begin{equation}\label{eqn-v1-shift}
v_{1}= (x_{p+1}+a_{c_{p}}) + (x_{p+2}+a_{c_{p+1}}) + \cdots + (x_{q-1}+a_{c_{q-2}}) + (x_p+y_q) \,.
\end{equation}
In our example this involves rewriting the original cluster (\ref{eqn-cluster-ex-short}) in the cyclically shifted form
\begin{equation}\label{eqn-cluster-ex-short-shifted}
3_1 3_2 4_4 6_7 6_8 6_9 (4_3\+5_5\+6_{6}\+7_{10}\+8_{11}\+\,\blue{3+9'}) 9_{13} 9_{14} 10'_{15} 10_{16}\,.
\end{equation}
where $s$ and $s'$ signify just $x_s$ and $y_s$, respectively.

In the general case, collecting together all the terms in $(x_p+y_q)$, that is summing $u_{k-1}(x_p+y_q)w_{m-k}$ over all possible 
$u_{k-1}$and $w_{m-k}$, yields
\begin{equation}
\resizebox{.9999 \textwidth}{!}
{ $
\begin{array}{l}
 (x_p+y_q)\ q_{k-1}(x_p, \sh x_{p+1}, \sh x_{p+2}, \cdots, \sh x_{q-1}|\a) \times 
         \ q_{m-k}(\sh^{q-p+1} x_q, y_{q+1}, x_{q+1}, \ldots, y_n, x_n|\a) \cr
						=(x_p+y_q)  q_{m-1} (x_p, \sh x_{p+1}, \cdots, \sh x_{q-1}, \sh x_q, y_{q+1}, x_{q+1}, \ldots, y_n, x_n|\a)\,, 
\end{array}
$  }						
\end{equation}
as required to produce the right hand side of (\ref{eqn-lem-b-proof}).

The remaining terms are of the form $u_{k-1}(x_r+a_{c_{r-1}})w_{m-k}$ for $r=p+1,p+2,\ldots,q-1$ 
with $(x_{r}+a_{c_{r-1}}) = (x_r+a_{r-p+\ell})$ for some $\ell$ determined by the requirement that $(x_s\+a_{s-p+\ell})$ 
is the rightmost factor of $u_{k-1}$ with $s<r$.
This ensures that if the rightmost factor $(x_s+a_{s-p+\ell})$ of $u_{k-1}$ is replaced by $(x_{r}+a_{c_{r-1}})$ to
give $u''_{k-1}$ then the result $u''_{k-1}$ is still a term in $q_{k-1}(x_p, \sh x_{p+1}, \sh x_{p+2}, \cdots, \sh x_{q-1}|\a)$
with $u''_{k-1}$ following $u_{k-1}$ in reverse lexicographic order, and $u''_{k-1}(x_s+a_{s-p+\ell})w_{m-k}$ constituting
a legitimate term in what we call a descendent cluster.

In our example (\ref{eqn-cluster-ex}), if we drop the factor $w_{m-k}$ which appears in every term as an innocent bystander, 
we have the following picture of ancestral clusters ($A_a$), the main cluster ($M$) and descendent clusters ($D_d$)
\begin{equation}
\resizebox{.99 \textwidth}{!}
{ $
\vcenter{\hbox{
\begin{tikzpicture}[x={(0in,-0.25in)},y={(0.75in,0in)}] 
\draw(0,0)node{$A_1~~ 3_1 3_2 3_3 6_7 6_8 6_9 (3_4\+4_5\+5_6\+6_{10}\+7_{11}\+8_{12}\+9'_{12})
=3_1 3_2 \red{3_3} 6_7 6_8 6_9 (\red{4_4}\+5_5\+6_{6}\+7_{10}\+8_{11}\+\blue{3\+9'})$};
\draw(2,0)node{$A_2~~3_1 3_2 4_4 4_5 6_8 6_9 (3_3\+4_6\+5_7\+6_{10}\+7_{11}\+8_{12}\+9'_{12})
=3_1 3_2 4_4 \red{4_5} 6_8 6_9 (4_3\+5_6\+\red{6_{7}}\+7_{10}\+8_{11}\+\blue{3\+9'})$};
\draw(4,0)node{$A_3~~3_1 3_2 4_4 5_6 6_8 6_9 (3_3\+4_5\+5_7\+6_{10}\+7_{11}\+8_{12}\+9'_{12})
=3_1 3_2 4_4 \red{5_6} 6_8 6_9 (4_3\+5_5\+\red{6_{7}}\+7_{10}\+8_{11}\+\blue{3\+9'})$};
\draw(7,0)node{$M~~3_1 3_2 4_4 6_7 6_8 6_9 (\red{3_3\+4_5\+5_6}\+\blue{6_{10}\+7_{11}\+8_{12}\+9'_{12}})
=3_1 3_2 4_4 6_7 6_8 6_9 (\red{4_3\+5_5\+6_{6}\+7_{10}\+8_{11}}\+\blue{3\+9'})$};
\draw(10,0)node{$D_1~~3_1 \red{4_3} 4_4 6_7 6_8 6_9 (\red{3_2}\+4_5\+5_6\+6_{10}\+7_{11}\+8_{12}\+9'_{12})
=3_1 4_3 4_4 6_7 6_8 6_9 (4_2\+5_5\+6_{6}\+7_{10}\+8_{11}\+\blue{3\+9'})$};
\draw(12,0)node{$D_2~~3_1 3_2 \red{5_5} 6_7 6_8 6_9 (3_3\+\red{4_4}\+5_6\+6_{10}\+7_{11}\+8_{12}\+9'_{12})
=3_1 3_2 5_5 6_7 6_8 6_9 (4_3\+5_4\+6_{6}\+7_{10}\+8_{11}\+\blue{3\+9'})$};
\draw(14,0)node{$D_3~~3_1 3_2 \red{6_6} 6_7 6_8 6_9 (3_3\+\red{4_4}\+5_6\+6_{10}\+7_{11}\+8_{12}\+9'_{12})
=3_1 3_2 6_6 6_7 6_8 6_9 (4_3\+5_4\+6_{6}\+7_{10}\+8_{11}\+\blue{3\+9'})$};
\draw(16,0)node{$D_4~~3_1 3_2 4_4 6_7 6_8 \red{7_{10}} (3_3\+4_5\+5_6\+\red{6_{9}}\+7_{11}\+8_{12}\+9'_{12})
=3_1 3_2 4_4 6_7 6_8 7_{10} (4_3\+5_5\+6_{6}\+7_{9}\+8_{11}\+\blue{3\+9'})$};
\draw(18,0)node{$D_5~~3_1 3_2 4_4 6_7 6_8 \red{8_{11}} (3_3\+4_5\+5_6\+\red{6_{9}}\+7_{10}\+8_{12}\+9'_{12})
=3_1 3_2 4_4 6_7 6_8 8_{11}(4_3\+5_5\+6_{6}\+7_{9}\+8_{10}\+\blue{3\+9'})$};
\draw[-latex](0+0.4,1.75)--(7-0.4,-2.10);
\draw[-latex](2+0.4,2.35)--(7-0.4,-1.75);
\draw[-latex](4+0.4,2.35)--(7-0.4,-1.45);
\draw[-latex](7+0.4,1.75)--(10-0.4,-2.10);
\draw[-latex](7+0.4,2.05)--(12-0.4,-1.75);
\draw[-latex](7+0.4,2.35)--(14-0.4,-1.75);
\draw[-latex](7+0.4,2.65)--(16-0.4,-1.05);
\draw[-latex](7+0.4,2.95)--(18-0.4,-1.05);
\end{tikzpicture}
}}
\\ \nonumber
$ 
}
\end{equation}

This process may be iterated, as in the following example for the case $n\geq 4$, $p=1$, $q=4$ and $m=4$. In this case
there are $4$ trees; one for each $k=1,2,3$ and $4$. For each $k$ the clusters are all of the of the form
$u_{k-1}v_1w_{m-k}$ with $u_{k-1}$ a term in $q_{k-1}(x_1,\sh x_2,\sh x_3)$,
$v_1=(x_1+a_{c_1}+x_2+a_{c_2}+x_3+a_{c_3}+y_4-a_{c_3})$ and $w_{4-k}$ a term in $q_{4-k}(\sh^2 x_4,y_5,x_5,\ldots,y_n,x_n)$. 
These may be drawn up for each $k$ in reverse lexicographic order as shown below in the columns on the left. Once again we have adopted 
the notation $s_t=(x_s+a_t)$ and $s'_t=(y_s-a_t)$ with $s=x_s$ and $s'=y_s$ and we have omitted the common factors of $w_{m-k}$. 
Then for each cluster the cyclic shift of subscripts on the $a$'s has been carried out, deleting the pairs $a_{k+2}$ and $-a_{k+2}$ and separating 
off the terms in $(x_1+y_4)$ to give the contributions on the right.

\begin{equation}
\vcenter{\hbox{
\begin{tikzpicture}[x={(0in,-0.20in)},y={(0.75in,0in)}] 
\draw(0,-1.75)node{$k=1$:};
 \draw(1,1)node{$(\blue{1_1+2_2+3_3+4'_3})\ =\ \red{(2_1+3_2)}\ +\ \blue{(1+4')}$}; 
\end{tikzpicture}
}}
\\ \nonumber
\end{equation}

\begin{equation}
\vcenter{\hbox{
\begin{tikzpicture}[x={(0in,-0.20in)},y={(0.75in,0in)}] 
\draw(3,-1.75)node{$k=2$:};
 \draw(4,1)node{$1_1 (\blue{1_2+2_3+3_4+4'_4})\ =\ 1_1 (2_2+3_3)\ +\ \blue{1_1 (1+4')}$}; 
 \draw(6,1)node{$2_2 (\red{1_1}+\blue{2_3+3_4+4'_4})\ =\ 2_2 (2_1+3_3)\ +\ \blue{2_2 (1+4')}$}; 
 \draw(8,1.7)node{$\red{2_1 2_2}$};
 \draw(10,1)node{$3_3 (\red{1_1+2_2}+\blue{3_4+4'_4})\ =\ 3_3 (2_1+3_2)\ +\ \blue{3_3 (1+4')}$}; 
 \draw(12,2.0)node{$\red{2_1 3_3 + 3_2 3_3}$}; 
\draw[-latex](4+0.5,1.35)--(6-0.5,-0.85); \draw[-latex](4+0.5,1.75)--(10-0.5,-0.85);
\draw[-latex](6+0.5,1.40)--(8-0.5,1.55); \draw[-latex](6+0.5,1.75)--(10-0.5,-0.40);
\draw[-latex](10+0.5,1.40)--(12-0.5,1.55); \draw[-latex](10+0.5,1.85)--(12-0.5,2.15);
\end{tikzpicture}
}}
\\ \nonumber
\end{equation}

\begin{equation}
\vcenter{\hbox{
\begin{tikzpicture}[x={(0in,-0.20in)},y={(0.75in,0in)}] 
\draw(14,-1.75)node{$k=3$:};
\draw(15,1)node{$1_1 1_2 (\blue{1_3+2_4+3_5+4'_5})\ =\ 1_1 1_2 (2_3+3_4)\ +\ \blue{1_1 1_2 (1+4')}$}; 
\draw(17,1)node{$1_1 2_3 (\red{1_2}+\blue{2_4+3_5+4'_5})\ =\ 1_1 2_3  (2_2+3_4)\ +\ \blue{1_1 2_3 (1+4')}$}; 
\draw(19,1)node{$2_2 2_3 (\red{1_1}+\blue{2_4+3_5+4'_5})\ =\ 2_2 2_3 (2_1+3_4)\ +\ \blue{2_2 2_3 (1+4')}$}; 
\draw(21,1.9)node{$\red{2_1 2_2 2_3}$}; 
\draw(23,1)node{$1_1 3_4 (\red{1_2+2_3}+\blue{3_5+4'_5})\ =\ 1_1 3_4 (2_2+3_3)\ +\ \blue{1_1 3_4 (1+4')}$}; 
\draw(25,1)node{$2_2 3_4 (\red{1_1+2_3}+\blue{3_5+4'_5})\ =\ 2_2 3_4 (2_1+3_3)\ +\ \blue{2_2 3_4 (1+4')}$}; 
\draw(27,1.9)node{$\red{2_1 2_2 3_4}$}; 
\draw(29,1)node{$3_3 3_4  (\red{1_1+2_2}+\blue{3_5+4'_5})\ =\ 3_3 3_4 (2_1+3_2)\ +\ \blue{3_3 3_4 (1+4')}$}; 
\draw(31,2.3)node{$\red{2_1 3_3 3_4 + 3_2 3_3 3_4 }$}; 
\draw[-latex](15+0.5,1.45)--(17-0.5,-0.95); \draw[-latex](15+0.5,1.90)--(23-0.5,-0.95);
\draw[-latex](17+0.5,1.45)--(19-0.5,-0.95); \draw[-latex](17+0.5,1.90)--(23-0.5,-0.50);
\draw[-latex](19+0.5,1.45)--(21-0.5,1.65); \draw[-latex](19+0.5,1.90)--(25-0.5,-0.50);
\draw[-latex](23+0.5,1.45)--(25-0.5,-0.95); \draw[-latex](23+0.5,1.90)--(29-0.5,-0.95);
\draw[-latex](25+0.5,1.45)--(27-0.5,1.65); \draw[-latex](25+0.5,1.90)--(29-0.5,-0.50);
\draw[-latex](29+0.5,1.45)--(31-0.5,1.65); \draw[-latex](29+0.5,1.90)--(31-0.5,2.45);
\end{tikzpicture}
}}
\\ \nonumber
\end{equation}

\begin{equation}
\vcenter{\hbox{
\begin{tikzpicture}[x={(0in,-0.25in)},y={(0.75in,0in)}] 
\draw(0,-1.75)node{$k=4$:};
 \draw(1,1)node{$1_1 1_2 1_3 (\blue{1_4+2_5+3_6+4'_6})\ =\ 1_1 1_2 1_3 (2_4+3_5)\ +\ \blue{1_1 1_2 1_3 (1+4')}$}; 
 \draw(3,1)node{$1_1 1_2 2_4 (\red{1_3}+\blue{2_5+3_6+4'_6})\ =\ 1_1 1_2 2_4 (2_3+3_5)\ +\ \blue{1_1 1_2 2_4 (1+4')}$}; 
 \draw(5,1)node{$1_1 2_3 2_4 (\red{1_2}+\blue{2_5+3_6+4'_6})\ =\ 1_1 2_3 2_4 (2_2+3_5)\ +\ \blue{1_1 2_3 2_4 (1+4')}$}; 
 \draw(7,1)node{$2_2 2_3 2_4 (\red{1_1}+\blue{2_5+3_6+4'_6})\ =\ 2_2 2_3 2_4 (2_1+3_5)\ +\ \blue{2_2 2_3 2_4 (1+4')}$}; 
 \draw(9,1.9)node{$\red{2_1 2_2 2_3 2_4}$}; 
\draw(11,1)node{$1_1 1_2 3_5 (\red{1_3+2_4}+\blue{3_6+4'_6})\ =\ 1_1 1_2 3_5 (2_3+3_4)\ +\ \blue{1_1 1_2 3_5 (1+4')}$}; 
\draw(13,1)node{$1_1 2_3 3_5 (\red{1_2+2_4}+\blue{3_6+4'_6})\ =\ 1_1 2_3 3_5 (2_2+3_4)\ +\ \blue{1_1 2_3 3_5 (1+4')}$}; 
\draw(15,1)node{$2_2 2_3 3_5 (\red{1_1+2_4}+\blue{3_6+4'_6})\ =\ 2_2 2_3 3_5 (2_1+3_4)\ +\ \blue{2_2 2_3 3_5 (1+4')}$}; 
\draw(17,1.9)node{$\red{2_1 2_2 2_3 3_5}$}; 
\draw(19,1)node{$1_1 3_4 3_5 (\red{1_2+2_3}+\blue{3_6+4'_6})\ =\ 1_1 3_4 3_5 (2_2+3_3)\ +\ \blue{1_1 3_4 3_5 (1+4')}$}; 
\draw(21,1)node{$2_2 3_4 3_5 (\red{1_1+2_3}+\blue{3_6+4'_6})\ =\ 2_2 3_4 3_5 (2_1+3_3)\ +\ \blue{2_2 3_4 3_5 (1+4')}$}; 
\draw(23,1.9)node{$\red{2_1 2_2 3_4 3_5}$}; 
\draw(25,1)node{$3_3 3_4 3_5 (\red{1_1+2_2}+\blue{3_6+4'_6})\ =\ 3_3 3_4 3_5 (2_1+3_2)\ +\ \blue{3_3 3_4 3_5 (1+4')}$}; 
\draw(27,2.0)node{$\red{2_1 3_3 3_4 3_5 + 3_2 3_3 3_4 3_5}$}; 
\foreach \i in{1,3,5,11,13,19} \draw[-latex](\i+0.4,1.55)--(\i+1.5,-1.05);
\draw[-latex](1+0.4,1.95)--(11-0.5,-1.05);
\draw[-latex](11+0.4,1.95)--(19-0.5,-1.05);
\draw[-latex](19+0.4,1.95)--(25-0.5,-1.05);
\foreach \i in{3,5,7} \draw[-latex](\i+0.4,1.95)--(\i+7.5,-0.65);
\foreach \i in{13,15} \draw[-latex](\i+0.4,1.95)--(\i+5.5,-0.65);
\foreach \i in{21} \draw[-latex](\i+0.4,1.95)--(\i+3.5,-0.65);
\foreach \i in{7,15,21} \draw[-latex](\i+0.4,1.55)--(\i+1.5,1.85);
\foreach \i in{25} \draw[-latex](\i+0.4,1.55)--(\i+1.5,1.45);
\foreach \i in{25} \draw[-latex](\i+0.4,1.95)--(\i+1.5,2.45);
\end{tikzpicture}
}}
\\ \nonumber
\end{equation}

For each $k$ the terms on the left in blue are those of $q_{k-1}(x_1,\sh x_2, \sh x_3, y_4|\a)$ 
and those in red are nonstandard and correspond, as indicated by the arrows, to the children of the terms on the right
with cyclically permuted shifts. 
On the right the terms in blue are those of $(x_1+y_4)q_{k-1}(x_1,\sh x_2,\sh x_3 |\a)$
and the remaining surviving terms on the right, known as ``leaves'' and shown in red, constitute $q_{k-1}(x_2,\sh x_3|\a)$. 
If one inserts the omitted common factors of $w_{m-k}$ and sums over these 
this tabulation serves to verify the validity of (\ref{eqn-lem-b-proof}) in the case $m=4$, $p=1$ and $q=4$ for any $n\geq 4$.

To complete the proof of our technical lemma, it remains to show that in the general case the ``leaves'', as exemplified above,
do indeed constitute $q_m(x_{p+1}, \sh x_{p+2}, \cdots,$ $ \sh x_{q-1}, \sh  x_q|\a)$.
To see this consider clusters of the form $u_{k-1}v_1w_{m-k}$ with 
\begin{equation}
      v_1= (x_p+a_1) + (x_{p+1}+a_2)  +\cdots + (x_r +a_r) + (x_{r+1}+ a_{c_{r+1}}) + \cdots (x_{q-1}+a_{c_{q-1}}) + (y_{q}-a_{c_{q-1}})\,,
\end{equation}
where $r$ is the maximum value of $t$ such that $(x_t+a_{c_t})=(x_t+a_t)$.
Such a cluster is valid if and only if $u_{k-1}$ consists of products of factors of the form $(x_i+a_j)$ with $i>r$.
The usual cyclic permutation of shifts and cancellation gives
\begin{equation}
      v_1=  (x_{p+1}+a_1) + (x_{p+2}+a_2)  + \cdots + (x_r+a_{r-1}) + (x_{r+1}+ a_{r}) + \cdots (x_{q-1}+a_{c_{q-2}}) + (x_p + y_{q})\,.
\end{equation}
Combining the first $r-p+1$ of these terms with each $u_{k-1}$ of the required form gives $r-p+1$ leaves belonging to
$q_{k}(x_{p+1},\sh x_{p+2}, \ldots, \sh x_{q-1}, y_q)|\a)$. Summing these leaves over all $u_{k-1}$ of the required form
gives $q_{k}(x_{p+1}, \sh x_{p+2}, \ldots, \sh x_{q-1}, y_q|\a)$. Each of the remaining terms $(x_t+a_{c_{t-1}})$ in $v_1$
gives rise in the usual way to a child in a lower cluster since the fact that $c_{t-1}>t-1$ ensures that an 
exchange is possible with a factor in $u_{k-1}$.

Finally, inserting the omitted factors $w_{m-k}$ and summing the leaves over $k=1,2,\ldots,n$ yields 
$q_m(x_{p+1},\sh x_{p+2}, \ldots, \sh x_{q-1}, y_q, \sh x_q, y_{q+1},x_{q+1},\ldots,y_n,x_n)|\a)$,
as is required to complete the proof of (\ref{eqn-lem-b-proof}) and consequently that of Lemma~\ref{lem-technical}
and Theorem~\ref{thm-main}.
\qed
\bigskip


\section{Corollaries}
\label{sec:corollaries}

First we note a corollary that is easily described, namely Lemma 4.10 of Ikeda {\em et al.}~\cite{ikeda} with the 
parameters $\a$ added rather than subtracted. 
\begin{Corollary}
Let $\mu$ be a partition of length $\ell(\mu)\leq n$ and $\delta= (n, n-1, \ldots, 1)$, so that
$\lambda=\mu+\delta$ is a strict partition of length $\ell(\lambda)=n$. Then
for $\a=(a_1,a_2,\ldots)$ we have:   
\begin{equation}
Q_\lambda(\x|\a) = \prod_{1\leq i \leq n} 2x_i \ \prod_{1\leq i<j \leq n} (x_i + x_j)\ s_{\mu} (\x|\a).
\end{equation}
\end{Corollary}
\noindent{\bf Proof}: One merely sets $\y=\x$ in equation (\ref{eqn-tok-fact-Qsfn}) of Theorem \ref{thm-main}. \qed

To make contact with other results it is necessary to introduce and relate a number of combinatorial constructs 
that are all in bijective correspondence with unprimed shifted tableaux, (USTx), namely 
strict Gelfand-Tsetlin patterns, (GTPs), certain alternating sign matrices, (ASMs), 
compass point matrices, (CPMs) and square-ice configurations, (SICs).


A Gelfand-Tsetlin pattern $G$ of size $n$ is a triangular array of non-negative integers $m_{ij}$ of the form
\begin{equation}\label{eqn-GT}
G = \left(
\begin{array}{ccccccccccccc}
 m_{n1}&&m_{n2}&&\cdots&&m_{n,n-1}&&m_{nn}\cr
&\ddots&&\ddots&&\ddots&&\cdots\cr
&&m_{31}&&m_{32}&&m_{33}\cr
&&&m_{21}&&m_{22}\cr
&&&&m_{11}\cr
\end{array}
\right)
\end{equation}
subject to the {\it betweenness conditions}
\begin{equation}\label{eqn-between}
    m_{i,j} \geq m_{i-1,j} \geq m_{i,j+1} \quad\hbox{for $i=2,3,\ldots,n$ and $j=1,2,\ldots,i-1$}\,.
\end{equation}
It follows that each row is a partition.
A Gelfand-Tsetlin pattern is said to be {\it strict} if
\begin{equation}\label{eqn-strict}
    m_{ij}>m_{i,j+1} \quad\hbox{for $i=1,2,\ldots,n$ and $j=1,2,\ldots,n-i$}\, ,
\end{equation}
in which case each row is a strict partition. 
A strict Gelfand-Tsetlin is sometimes called a {\it monotone triangle}~\cite{Okadapartial}. 

For each strict partition $\lambda$ of length $\ell(\lambda)=n$ let ${\cal G}^\lambda(n)$ be the set of all 
strict GTPs 
$G$ with top row $\lambda$, that is $m_{ni}=\lambda_i$ for
$i=1,2,\ldots,n$. These are in bijective correspondence with all USTx 
$S\in{\cal S}^\lambda[\n]$,
where the correspondence is defined by
\begin{equation}\label{eqn-mij}
  m_{ij} = \hbox{number of entries $\leq i$ in row $j$ of $S$} \,. 
\end{equation}
Conversely,
\begin{equation}
  s_{ij} = \left\{ 
	     \begin{array}{cl}
			     1 &\hbox{if $i=1$ and $j\leq m_{11}$}\cr
					 k &\hbox{if $i>1$ and $m_{k-1,i}<j\leq m_{ki}$ for each $k=2,\ldots,n$}\cr
				\end{array} \right.
\end{equation}
It is straightforward to check that with the constraints (\ref{eqn-between}) and (\ref{eqn-strict}) the 
conditions S1-S3 of section~\ref{sec:background} are automatically satisfied and {\it vice versa}. 

Next we turn to 
ASMs. For each strict partition $\lambda$ of length $\ell(\lambda)=n$ and breadth 
$\lambda_1=m$ let ${\cal A}^\lambda$ be the set of all $n\times m$ matrices $A=(a_{ij})_{1\leq i\leq n,1\leq j\leq m}$
with $a_{i,j}\in\{1,0,-1\}$ such that
\begin{itemize}
\item[A1] the non-zero entries alternate in sign across each row and down each column;
\item[A2] the rightmost non-zero entry in each row is $1$;
\item[A3] the topmost non-zero entry in any column is $1$;
\item[A4] $\sum_{j=1}^m a_{ij}=1$ for $i=1,2,\ldots,n$;
\item[A5] $\sum_{i=1}^n a_{ij}=1$ if $j=\lambda_k$ for some $k$ and $0$ otherwise. 
\end{itemize}

These too are in bijective correspondence with strict GTPs 
$G\in{\cal G}^\lambda$ with the correspondence defined by~\cite{MRR}
\begin{equation}
   a_{ij}=\left\{  
	           \begin{array}{cl}
						  1&\hbox{if $i=1$ and the $1$st row of $G$ contains $j$;}\cr
						  1&\hbox{if $i>1$ and the $i$th row of $G$ contains $j$ but the $(i-1)$th does not;}\cr
							-1&\hbox{if $i>1$ and the $(i-1)$th row of $G$ contains $j$ but the $i$th does not;}\cr
							0&\hbox{otherwise,} \cr
							\end{array} \right.
\end{equation}
where it might be noted that the rows of $G$ are counted from bottom to top and those of $A$ from top to bottom. Similarly, the bijective
correspondence with USTx 
$S\in{\cal S}^\lambda[\n]$ is defined by
\begin{equation}
   a_{ij}=\left\{  
	           \begin{array}{cl}
						  1&\hbox{if $j=m$ and the $m$th diagonal of $S$ contains $i$;}\cr
						  1&\hbox{if $j<m$ and the $j$th diagonal of $S$ contains $i$ but  $(j+1)$th does not;}\cr
							-1&\hbox{if $j<m$ and the $(j+1)$th diagonal of $S$ contains $i$ but  $j$th does not;}\cr
							0&\hbox{otherwise,} \cr
							\end{array} \right.
\end{equation}

As emphasised elsewhere~\cite{HKbij}, to each ASM we can associate both 
a CPM 
and an SIC 
of the $6$-vertex model.
We define the CPMs $C\in{\cal C}^\lambda$ corresponding to $A\in{\cal A}^\lambda$ to be those matrices
obtained by mapping the entries $1$ and $-1$ in $A$ to $\WE$ and $\NS$, respectively, 
and mapping an entry $0$ in $A$ to one or other of $\NE$, $\SE$, $\NW$ or $\SW$ 
in accordance with the compass point arrangements of the nearest non-zero  
neighbours of the $0$, as specified in the tabulation (\ref{eqn-ASM-CPM-SIC}).

Each SIC takes the form of a planar grid consisting of vertices and directed edges.  
Each vertex has four edges, two incoming and two outgoing, resulting in six vertex 
configurations that may be constructed from the six possible entries $\XY$ 
of a CPM 
by attaching to a vertex two incoming edges from the directions $\X$ and $\Y$
with the other two edges outgoing, as shown in the fourth row of table (\ref{eqn-ASM-CPM-SIC}).

\begin{equation}\label{eqn-ASM-CPM-SIC}
\begin{array}{|c|c|c|c|c|c|c|}
\hline
&&&&&&\cr
\hbox{ASM} & 1 &\ov1 & 0 & 0 &0 &0 \cr
\hline
&&&&&&\cr
                &\begin{array}{c}\ov1\cr \ov1~~{\bf1}~~\ov1\cr\ov1\end{array}
								&\begin{array}{c} 1\cr 1~~{\bf\ov1}~~1\cr 1\end{array}
								&\begin{array}{c} 1\cr 1~~{\bf0}~~\ov1\cr\ov1\end{array}
								&\begin{array}{c}\ov1\cr 1~~{\bf0}~~\ov1\cr 1\end{array}
								&\begin{array}{c} 1\cr \ov1~~{\bf0}~~1\cr\ov1\end{array}
								&\begin{array}{c}\ov1\cr \ov1~~{\bf0}~~1\cr 1\end{array}\cr
&&&&&&\cr								
\hline
\hbox{CPM} & \WE &\NS & \NE & \SE & \NW & \SW \cr
\hline
&&&&&&\cr								

\hbox{SIC} & 
\vcenter{\hbox{
\begin{tikzpicture}[x={(0in,-0.25in)},y={(0.25in,0in)}] 
\draw(1,2)node{$\bullet$};
\draw[thick,->] (1,1)to(1,1.6); \draw[thick] (1,1.6)to(1,2);
\draw[thick] (1,2)to(1,2.4); \draw[thick,<-] (1,2.4)to(1,3);
\draw[thick](0,2)to(0.4,2); \draw[thick,<-] (0.4,2)to(1,2);
\draw[thick,->](1,2)to(1.6,2); \draw[thick] (1.6,2)to(2,2);
\end{tikzpicture}
}}
& 
\vcenter{\hbox{
\begin{tikzpicture}[x={(0in,-0.25in)},y={(0.25in,0in)}] 
\draw(1,2)node{$\bullet$};
\draw[thick] (1,1)to(1,1.4); \draw[thick,<-] (1,1.4)to(1,2);
\draw[thick,->] (1,2)to(1,2.6); \draw[thick] (1,2.6)to(1,3);
\draw[thick,->](0,2)to(0.6,2); \draw[thick] (0.6,2)to(1,2);
\draw[thick](1,2)to(1.4,2); \draw[thick,<-] (1.4,2)to(2,2);
\end{tikzpicture}
}}
& 
\vcenter{\hbox{
\begin{tikzpicture}[x={(0in,-0.25in)},y={(0.25in,0in)}] 
\draw(1,2)node{$\bullet$};
\draw[thick] (1,1)to(1,1.4); \draw[thick,<-] (1,1.4)to(1,2); 
\draw[thick] (1,2)to(1,2.4); \draw[thick,<-] (1,2.4)to(1,3); 
\draw[thick,->](0,2)to(0.6,2); \draw[thick] (0.6,2)to(1,2);  
\draw[thick,->](1,2)to(1.6,2); \draw[thick] (1.6,2)to(2,2);  
\end{tikzpicture}
}}
&\vcenter{\hbox{
\begin{tikzpicture}[x={(0in,-0.25in)},y={(0.25in,0in)}] 
\draw(1,2)node{$\bullet$}; 
\draw[thick] (1,1)to(1,1.4); \draw[thick,<-] (1,1.4)to(1,2);
\draw[thick] (1,2)to(1,2.4); \draw[thick,<-] (1,2.4)to(1,3);
\draw[thick](0,2)to(0.4,2); \draw[thick,<-] (0.4,2)to(1,2);
\draw[thick](1,2)to(1.4,2); \draw[thick,<-] (1.4,2)to(2,2);
\end{tikzpicture}
}}
&\vcenter{\hbox{
\begin{tikzpicture}[x={(0in,-0.25in)},y={(0.25in,0in)}] 
\draw(1,2)node{$\bullet$}; 
\draw[thick,->] (1,1)to(1,1.6); \draw[thick] (1,1.6)to(1,2);
\draw[thick,->] (1,2)to(1,2.6); \draw[thick] (1,2.6)to(1,3);
\draw[thick,->](0,2)to(0.6,2); \draw[thick] (0.6,2)to(1,2);
\draw[thick,->](1,2)to(1.6,2); \draw[thick] (1.6,2)to(2,2);
\end{tikzpicture}
}}
&\vcenter{\hbox{
\begin{tikzpicture}[x={(0in,-0.25in)},y={(0.25in,0in)}] 
\draw(1,2)node{$\bullet$};
\draw[thick,->] (1,1)to(1,1.6); \draw[thick] (1,1.6)to(1,2);
\draw[thick,->] (1,2)to(1,2.6); \draw[thick] (1,2.6)to(1,3);
\draw[thick](0,2)to(0.4,2); \draw[thick,<-] (0.4,2)to(1,2);
\draw[thick](1,2)to(1.4,2); \draw[thick,<-] (1.4,2)to(2,2);
\end{tikzpicture}
}}
\cr
&&&&&&\cr								
\hline
\end{array}
\end{equation}

In this table, and in the example that follows in (\ref{eqn-ASCGI}), each symbol $\ov{1}$ is to be interpreted as an ASM entry $-1$. The first 
row of the table specifies an ASM entry that is further characterised in the second row by its four outer $1$s and $\ov{1}$s 
indicating the values of the nearest non-zero ASM entries to its north, east, south or west, where any missing non-zero neighbour to 
the east or north is taken to be $\ov{1}$.

The admissible square ice configurations $I\in{\cal I}^\lambda$ 
are defined to be those constructed from the six vertices in the form of $n\times m$ grids, with $n=\ell(\lambda)$, $m=\lambda_1$,
for which the boundary horizontal edges are all incoming and the boundary vertical edges
are all outgoing except for those on the lower boundary in columns that do not correspond to a part of $\lambda$.
The maps defined by (\ref{eqn-ASM-CPM-SIC}) from ASMs to CPMs to SICs are easily shown to ensure that the $I\in{\cal I}^\lambda$ are 
in bijective correspondence with the $A\in{\cal A}^\lambda$.

The bijections between ${\cal A}^\lambda$, ${\cal S}^\lambda$,  ${\cal G}^\lambda$  and ${\cal I}^\lambda$ are all exemplified in (\ref{eqn-ASCGI}).

\begin{center}
\begin{equation}\label{eqn-ASCGI}
\begin{array}{l}
A =
\left[
\begin{array}{cccccc}
0&1&0&0&0&0\cr
1&\ov1&0&1&0&0\cr
0&0&0&0&1&0\cr
0&0&1&0&\ov1&1\cr
\end{array}
\right]
\Leftrightarrow
\left[
\begin{array}{cccccc}
1&1&0&0&0&0\cr
1&0&1&1&0&0\cr
1&1&1&1&1&0\cr
1&1&1&0&0&1\cr
\end{array}
\right]
\Leftrightarrow
\SYT{0.2in}{}{
 {1,1,2,2,3,4},
   {2,3,3,3},
     {3,4,4},
       {4}
}
= S
\cr\cr
{\hskip7em}\Updownarrow{\hskip5em}\searrow{\hskip-0.9em}\nwarrow
\cr\cr
{\phantom{A=}}
\left[
\begin{array}{cccccc}
0&1&0&0&0&0\cr
1&0&0&1&0&0\cr
1&0&0&1&1&0\cr
1&0&1&1&0&1\cr
\end{array}
\right]
{\hskip1em}
~~~~
C =
\left[
\begin{array}{cccccc}
\SW&\WE&\SE&\SE&\SE&\SE\cr
\WE&\NS&\SW&\WE&\SE&\SE\cr
\NW&\SW&\SW&\NW&\WE&\SE\cr
\NW&\SW&\WE&\NE&\NS&\WE\cr
\end{array}
\right]
\cr\cr
{\hskip7em}\Updownarrow{\hskip14em}\searrow{\hskip-0.9em}\nwarrow
\cr\cr
G = 
\left(\!
\begin{array}{ccccccccccccc}
6&&4&&3&&1\cr
&5&&4&&1\cr
&&4&&1\cr
&&&2\cr
\end{array}
\!\right)
{\hskip9em}
I = 
\vcenter{\hbox{
\begin{tikzpicture}[x={(0in,-0.2in)},y={(0.2in,0in)}] 
\foreach \i in {1,...,4}
\foreach \j in {1,...,6}
\draw(\i,\j)node{$\bullet$};
%
%
\draw[thick, ->] (1,0)to(1,0.6); \draw[thick] (1,0.6)to(1,1);
\draw[thick, ->] (1,1)to(1,1.6); \draw[thick] (1,1.6)to(1,2);
\draw[thick] (1,2)to(1,2.4); \draw[thick,<-] (1,2.4)to(1,3);
\draw[thick] (1,3)to(1,3.4); \draw[thick,<-] (1,3.4)to(1,4);
\draw[thick] (1,4)to(1,4.4); \draw[thick,<-] (1,4.4)to(1,5);
\draw[thick] (1,5)to(1,5.4); \draw[thick,<-] (1,5.4)to(1,6);
\draw[thick] (1,6)to(1,6.4); \draw[thick,<-] (1,6.4)to(1,7);
\draw[thick, ->] (2,0)to(2,0.6); \draw[thick] (2,0.6)to(2,1);
\draw[thick] (2,1)to(2,1.4);  \draw[ thick,<-] (2,1.4)to(2,2);
\draw[thick,->] (2,2)to(2,2.6);  \draw[thick] (2,2.6)to(2,3);
\draw[thick,->] (2,3)to(2,3.6);  \draw[ thick] (2,3.6)to(2,4);
\draw[thick] (2,4)to(2,4.4); \draw[thick,<-] (2,4.4)to(2,5);
\draw[thick] (2,5)to(2,5.4);  \draw[ thick,<-] (2,5.4)to(2,6);
\draw[thick] (2,6)to(2,6.4);  \draw[thick,<-] (2,6.4)to(2,7);
\draw[thick, ->] (3,0)to(3,0.6); \draw[thick] (3,0.6)to(3,1);
\draw[thick, ->] (3,1)to(3,1.6);  \draw[thick] (3,1.6)to(3,2);
\draw[thick, ->] (3,2)to(3,2.6);  \draw[thick] (3,2.6)to(3,3);
\draw[thick, ->] (3,3)to(3,3.6);  \draw[thick] (3,3.6)to(3,4);
\draw[thick, ->] (3,4)to(3,4.6);  \draw[thick] (3,4.6)to(3,5);
\draw[thick] (3,5)to(3,5.4);  \draw[thick,<-] (3,5.4)to(3,6);
\draw[thick] (3,6)to(3,6.4);  \draw[ thick,<-] (3,6.4)to(3,7);
\draw[thick, ->] (4,0)to(4,0.6); \draw[thick] (4,0.6)to(4,1);
\draw[thick,->] (4,1)to(4,1.6);  \draw[ thick] (4,1.6)to(4,2);
\draw[thick,->] (4,2)to(4,2.6);  \draw[ thick] (4,2.6)to(4,3);
\draw[thick] (4,3)to(4,3.4);  \draw[ thick,<-] (4,3.4)to(4,4);
\draw[thick] (4,4)to(4,4.4); \draw[thick,<-] (4,4.4)to(4,5);
\draw[thick,->] (4,5)to(4,5.6);  \draw[ thick] (4,5.6)to(4,6);
\draw[thick] (4,6)to(4,6.4);  \draw[ thick,<-] (4,6.4)to(4,7);
%
%
\draw[thick](0,1)to(0.4,1); \draw[thick,<-] (0.4,1)to(1,1); 
\draw[thick](0,2)to(0.4,2); \draw[thick,<-] (0.4,2)to(1,2);
\draw[thick](0,3)to(0.4,3); \draw[thick,<-] (0.4,3)to(1,3);
\draw[thick](0,4)to(0.4,4); \draw[thick,<-] (0.4,4)to(1,4);
\draw[thick](0,5)to(0.4,5); \draw[thick,<-] (0.4,5)to(1,5);
\draw[thick](0,6)to(0.4,6); \draw[thick,<-] (0.4,6)to(1,6);
\draw[thick](1,1)to(1.4,1); \draw[thick,<-] (1.4,1)to(2,1);
\draw[thick,->](1,2)to(1.6,2); \draw[thick] (1.6,2)to(2,2);
\draw[thick](1,3)to(1.4,3); \draw[thick,<-] (1.4,3)to(2,3);
\draw[thick](1,4)to(1.4,4); \draw[thick,<-] (1.4,4)to(2,4);
\draw[thick](1,5)to(1.4,5); \draw[thick,<-] (1.4,5)to(2,5);
\draw[thick](1,6)to(1.4,6); \draw[thick,<-] (1.4,6)to(2,6);
\draw[thick,->](2,1)to(2.6,1); \draw[thick] (2.6,1)to(3,1);
\draw[thick](2,2)to(2.4,2); \draw[thick,<-] (2.4,2)to(3,2);
\draw[thick](2,3)to(2.4,3); \draw[thick,<-] (2.4,3)to(3,3);
\draw[thick,->](2,4)to(2.6,4); \draw[thick] (2.6,4)to(3,4);
\draw[thick](2,5)to(2.4,5); \draw[thick,<-] (2.4,5)to(3,5);
\draw[thick](2,6)to(2.4,6); \draw[thick,<-] (2.4,6)to(3,6);
\draw[thick,->](3,1)to(3.6,1); \draw[thick] (3.6,1)to(4,1);
\draw[thick](3,2)to(3.4,2); \draw[thick,<-] (3.4,2)to(4,2);
\draw[thick](3,3)to(3.4,3); \draw[thick,<-] (3.4,3)to(4,3);
\draw[thick](3,4)to(3.4,4); \draw[thick,<-] (3.4,4)to(4,4);
\draw[thick,->](3,5)to(3.6,5); \draw[thick] (3.6,5)to(4,5);
\draw[thick](3,6)to(3.4,6); \draw[thick,<-] (3.4,6)to(4,6);
\draw[thick,->](4,1)to(4.6,1); \draw[thick] (4.6,1)to(5,1);
\draw[thick](4,2)to(4.4,2); \draw[thick,<-] (4.4,2)to(5,2);
\draw[thick,->](4,3)to(4.6,3); \draw[thick] (4.6,3)to(5,3);
\draw[thick,->](4,4)to(4.6,4); \draw[thick] (4.6,4)to(5,4);
\draw[thick](4,5)to(4.4,5); \draw[thick,<-] (4.4,5)to(5,5);
\draw[thick,->](4,6)to(4.6,6); \draw[thick] (4.6,6)to(5,6);
\end{tikzpicture}
}}
\end{array}
\end{equation}
\end{center}

As illustrated in the top row of (\ref{eqn-ASCGI}), the map from $A\in{\cal A}^\lambda$ to $S\in{\cal S}^\lambda[\n]$ may be constructed by first drawing up a cumulative 
row sum matrix of $1$s and $0$s by summing the entries of $A$ from right to left across its rows, and then filling the $j$th diagonal of $S$ from 
top-left to bottom-right with the row numbers of the $1$s appearing from top to bottom in the $j$th column of the cumulative row sum matrix. 
Similarly the bijective map from $A\in{\cal A}^\lambda$ to $G\in{\cal G}^\lambda$ may be constructed by first drawing up a cumulative column sum matrix $cs(A)$ of $1$s and $0$s by summing entries the entries of $A$ from top to bottom down its columns, and then filling the $i$th row of $G$ from left to right with the column numbers of the 
$1$s appearing from left to right in the $i$th row of the cumulative column sum matrix. This is illustrated in left hand column of (\ref{eqn-ASCGI}). 
Finally the map from $A\in{\cal A}^\lambda$ to $I\in{\cal I}^\lambda$ proceeds, as shown on 
the diagonal of (\ref{eqn-ASCGI}), by way of the compass point matrix $C$ in accordance with the six-vertex tabulation
of (\ref{eqn-ASM-CPM-SIC}). The simplicity of these maps makes it easy to check that they are all bijections.

In order to establish corollaries of our main result Theorem~\ref{thm-main} within the context of the above combinatorial objects
it is merely necessary to replace the sum over $P\in{\cal P}^\lambda(\n,\n')$ by sums over $K\in{\cal K}^\lambda$
with an appropriate identification of $\wgt(K)$ in the three cases ${\cal K}^\lambda={\cal A}^\lambda$, ${\cal G}^\lambda$, and
${\cal I}^\lambda$. 

The simplest case is that of ${\cal A}^\lambda$ for which :

\begin{Corollary}\label{cor-A}
Let $\lambda$ be a strict partition of length $\ell(\lambda)=n$ and breadth $\lambda_1=m$ and let $\a=(a_1,a_2,\ldots)$. 
Then for each $n\times m$ ASM $A\in{\cal A}^\lambda$ let $C(A)=(c_{ij})$ be the corresponding CPM. Then
\begin{equation}\label{eqn-tok-fact-A}
     \sum_{A\in{\cal A}^\lambda} \wgt(A)  = \prod_{i=1}^n x_i \ \prod_{1\leq i<j \leq n} (x_i + y_j) \ \ s_{\mu} (\x|\a)\,,
\end{equation}
where $\mu=\lambda-\delta$ with $\delta=(n,n-1,\ldots,1)$ and
\begin{equation}\label{eqn-wgtA}
      \wgt(A) = \prod_{i=1}^n x_i\  \prod_{i=1}^n\prod_{j=1}^m \wgt(c_{ij})
\end{equation}
with
\begin{equation}\label{eqn-wgtC}
\begin{array}{|c|c|c|c|c|c|c|c|}
\hline
\hbox{Entry} &&&&&& \cr
\hbox{at $(i,j)$} &\WE &\NS & \NE & \SE & \NW & \SW \cr
\hline
\wgt(c_{ij})& 1 & x_i+y_j &1 &1 & y_i-a_j &x_i+a_j\cr
\hline
\end{array}
\end{equation}							
\end{Corollary}
\bigskip

\noindent{\bf Proof}:~~
The right hand side of (\ref{eqn-tok-fact-A}) coincides with that of (\ref{eqn-tok-fact-PT}) so that all we have to show is that
\begin{equation}\label{eqn-AtoPtoS}
     \sum_{A\in{\cal A}^\lambda} \wgt(A)  =  \sum_{P\in{\cal P}^\lambda(\n,\n')} \wgt(P) = \sum_{S\in{\cal S}^\lambda(\n)} \wgt(S)
\end{equation}
where $\wgt(P)$ is defined by the left hand parts of (\ref{eqn-Pwgt-Twgt}) and (\ref{eqn-pwgt-twgt}). However the one-to-many map from $S\in{\cal S}^\lambda(\n)$
to $P\in{\cal P}^\lambda(\n,\n')$ is such that 
\begin{equation}
\sum_{P\in{\cal P}^\lambda(\n,\n')} \wgt(P) = \sum_{S\in{\cal S}^\lambda(\n)} \wgt(S)
\quad\hbox{where}\quad 
\wgt(S) = \sum_{(i,j)\in{\cal SF}^\lambda} \wgt(s_{ij})
\end{equation}
with
\begin{equation}\label{eqn-wgtS}
    \wgt(s_{ij})=\left\{   \begin{array}{ll} 
		                            x_k&\hbox{if $i=j$ and $s_{ii}=k$;}\cr
																x_k+a_{j-i}&\hbox{if $i<j$, $s_{ij}=k$ and $s_{i,j-1}=k$;}\cr
																y_k-a_{j-i}&\hbox{if $i<j$, $s_{ij}=k$ and $s_{i+1,j}=k$;}\cr
																x_k+y_k&\hbox{if $i<j$, $s_{ij}=k$, $s_{i,j-1}\neq k$ and $s_{i+1,j}\neq k$,}\cr
														\end{array}
		\right.
\end{equation} 
where the last case follows from the fact that $x_k+a_{j-i}+y_k-a_{j-i}=x_k+y_k$.
Now we only have to ensure that $\wgt(A)=\wgt(S)$ where $A$ and $S$ are related by the bijective map we have 
identified from $A\in{\cal A}^\lambda$ to $S\in{\cal S}^\lambda(\n)$.

The diagonal elements of any $S\in{\cal S}^\lambda(\n)$ are always $1,2,\ldots,n$ since they are strictly increasing down this diagonal.
Since in this case $i=j$ it follows from (\ref{eqn-wgtS}) that their contribution to $\wgt(S)$ is just the 
factor $x_1x_2\cdots x_n$ that appears on the right hand side of the expression (\ref{eqn-wgtA}) for $\wgt(A)$. 
To determine the remaining factors
it should be noted that the passage from $S$ to $C$ is such that each entry $k$ in diagonal $d>1$ is mapped to an entry 
$\SW$, $\NW$ or $\NS$ in row $k$ and column $d-1$ of $C$ according as there is another entry $k$ immediately to its left, another entry 
$k$ immediately below or no entry $k$ in diagonal $d-1$. It follows from (\ref{eqn-wgtS}) that the corresponding weights in $S$ are 
$x_k+a_{d-1}$, $y_k-a_{d-1}$ and $x_k+y_k$. Taking into account the shift from $d$ to $d-1$, and identifying $(k,d-1)$ with $(i,j)$ gives 
$\wgt(c_{ij})$ as tabulated in (\ref{eqn-wgtC}). Thus $\wgt(A)=\wgt(S)$ as required.
\qed
\bigskip

In our example (\ref{eqn-ASCGI}) this is illustrated by the following example in which $S$ and $C$ are shown on the left with their 
weights given by the product of all the entries on the right:
\begin{equation*}
\resizebox{.99 \textwidth}{!}
{ $
\begin{array}{l}
\SYT{0.2in}{}{
 {1,1,2,2,3,4},
   {2,3,3,3},
     {3,4,4},
       {4}
}
\!\mapsto\! 
\wideSYT{0.2in}{0.5in}{}{
 {x_1,x_1\!+\!a_1,x_2\!+\!y_2,x_2\!+\!a_3,y_3\!-\!a_4,x_4\!+\!y_4},
   {x_2,y_3\!-\!a_1,x_3\!+\!a_2,x_3\!+\!a_3},
     {x_3,y_4\!-\!a_1,x_4\!+\!a_2},
       {x_4},
}
\cr\cr
\left[
\begin{array}{cccccc}
\SW&\WE&\SE&\SE&\SE&\SE\cr
\WE&\NS&\SW&\WE&\SE&\SE\cr
\NW&\SW&\SW&\NW&\WE&\SE\cr
\NW&\SW&\WE&\NE&\NS&\WE\cr
\end{array}
\right]
\!\mapsto\!  
\wideYT{0.2in}{0.5in}{}{
 {x_1\!+\!a_1,1,1,1,1,1},
   {1,x_2\!+\!y_2,x_2+a_3,1,1,1},
     {y_3\!-\!a_1,x_3\!+\!a_2,x_3\!+\!a_3,y_3\!-\!a_4,1,1},
       {y_4\!-\!a_1,x_4\!+\!a_2,1,1,x_4\!+\!y_4,1},
}
\end{array}
 $ 
}
\end{equation*}

Thanks to the tabulation (\ref{eqn-ASM-CPM-SIC}) this corollary covers the cases ${\cal A}^\lambda$ and ${\cal I}^\lambda$.
It remains to consider the case ${\cal G}^\lambda$. 

\begin{Corollary}\label{cor-G}
Let $\lambda$ be a strict partition of length $\ell(\lambda)=n$ and let $\a=(a_0,a_1,a_2,\ldots)$ with $a_0=0$.
Then for each strict Gelfand-Tsetlin pattern $G\in{\cal G}^\lambda$, with entries $m_{ij}$ for 
$i=1,2,\ldots,n$ and $j=1,2,\ldots,i$,
\begin{equation}\label{eqn-tok-fact-G}
     \sum_{G\in{\cal G}^\lambda} \wgt(G)  = \prod_{i=1}^n x_i \ \prod_{1\leq i<j \leq n} (x_i + y_j) \ \ s_{\mu} (\x|\a)\,,
\end{equation}
where $\mu=\lambda-\delta$ with $\delta=(n,n-1,\ldots,1)$ and
\begin{equation}\label{eqn-wgtG}
\begin{array}{rcl}
 \wgt(G) &=&\ds \prod_{i=1}^n \prod_{k=0}^{m_{ii}-1} (x_i+a_k) \cr  
&&\ds	  ~\times~\prod_{i=2}^n\prod_{j=1}^{i-1} \big(\chi(B)(x_i+y_i)+\chi(R_k)(y_i-a_k)\big)\ \prod_{k=m_{i-1,j}+1}^{m_{ij}-1} (x_i+a_k)  \,,
\end{array}
\end{equation}
where $B:=m_{i,j}>m_{i-1,j}>m_{i,j+1}$, $R_k:=m_{i,j}>m_{i-1,j}=m_{i,j+1}=k$ and 
$\chi(P)$ is the truth function whereby $\chi(P)=1$ if the proposition $P$ is true, and $0$ otherwise.		
\end{Corollary}
\bigskip

\noindent{\bf Proof}:~~
As in the previous corollary, it is only necessary to establish that $\wgt(G)=\wgt(S)$ where $G$ and $S$ are related through 
the bijection between $G\in{\cal G}^\lambda$ and $S\in{\cal S}^\lambda(\n)$ that is defined by (\ref{eqn-mij}).
This states that the entry $m_{ij}$ in $G$ is the number of entries no greater than $i$ in row $j$ of the corresponding shifted tableau $S$. 
Thanks to the betweenness and strictness conditions
(\ref{eqn-between}) and (\ref{eqn-strict}) there are three cases to consider:
\begin{equation}\label{eqn-LRB}
\begin{array}{lll}
\hbox{(L)}&m_{i,j}=m_{i-1,j}>m_{i,j+1}\geq m_{i-1,j+1}&\vcenter{\hbox{
\begin{tikzpicture}[x={(0in,-0.2in)},y={(0.2in,0in)}] 
\draw[<-](-0.5,0)--(-0.5,4); \path(-0.4,4)--node{$m_{ij}$}(-0.4,6);  \draw[->](-0.5,6)--(-0.5,10);
\draw(0,0)--(0,10)--(1,10)--(1,0)--cycle;  \path(0.6,0)--node{$m_{i-1,j}=m_{ij}$}(0.6,10); 
\draw(1,1)--(1,5)--(2,5)--(2,1)--cycle; \path(1.6,1)--node{$m_{i-1,j+1}$}(1.6,5); \draw(1,5)--(1,7)--(2,7)--(2,5)--cycle; \path(1.5,5)--node{$i~i~i$}(1.5,7); 
\draw[<-](2.5,1)--(2.5,2.5); \path(2.6,3)--node{$m_{i,j+1}$}(2.6,5);  \draw[->](2.5,5.5)--(2.5,7);
\end{tikzpicture}
}}\cr
\hbox{(R)}&m_{i,j}>m_{i-1,j}=m_{i,j+1}>m_{i-1,j+1}&\vcenter{\hbox{
\begin{tikzpicture}[x={(0in,-0.2in)},y={(0.2in,0in)}] 
\draw[<-](-0.5,0)--(-0.5,4); \path(-0.4,4)--node{$m_{ij}$}(-0.4,6);  \draw[->](-0.5,6)--(-0.5,10);
\draw(0,0)--(0,7)--(1,7)--(1,0)--cycle;  \path(0.6,0)--node{$m_{i-1,j}$}(0.6,7); 
        \draw(0,7)--(0,8)--(1,8)--(1,7)--cycle; \path(0.5,7)--node{$i$}(0.5,8); 
				       \draw(0,8)--(0,10)--(1,10)--(1,8)--cycle; \path(0.5,8)--node{$i~i~i$}(0.5,10); 
\draw(1,1)--(1,5)--(2,5)--(2,1)--cycle; \path(1.6,1)--node{$m_{i-1,j+1}$}(1.6,5); 
         \draw(1,5)--(1,7)--(2,7)--(2,5)--cycle; \path(1.5,5)--node{$i~i~i$}(1.5,7); 
				       \draw(1,7)--(1,8)--(2,8)--(2,7)--cycle; \path(1.5,7)--node{$i$}(1.5,8); 
\draw[<-](2.5,1)--(2.5,1.7); \path(2.6,1.7)--node{$m_{i,j+1}=m_{i-1,j}$}(2.6,7.3);  \draw[->](2.5,7.3)--(2.5,8);
\end{tikzpicture}
}}\cr
\hbox{(B)}&m_{i,j}>m_{i-1,j}>m_{i,j+1}\geq m_{i-1,j+1}&\vcenter{\hbox{
\begin{tikzpicture}[x={(0in,-0.2in)},y={(0.2in,0in)}] 
\draw[<-](-0.5,0)--(-0.5,4); \path(-0.4,4)--node{$m_{ij}$}(-0.4,6);  \draw[->](-0.5,6)--(-0.5,10);
\draw(0,0)--(0,7)--(1,7)--(1,0)--cycle;  \path(0.6,0)--node{$m_{i-1,j}$}(0.6,7); 
        \draw(0,7)--(0,8)--(1,8)--(1,7)--cycle; \path(0.5,7)--node{$i$}(0.5,8); 
				       \draw(0,8)--(0,10)--(1,10)--(1,8)--cycle; \path(0.5,8)--node{$i~i~i$}(0.5,10); 
\draw(1,1)--(1,4)--(2,4)--(2,1)--cycle; \path(1.6,1)--node{$m_{i-1,j+1}$}(1.6,4); 
         \draw(1,4)--(1,6)--(2,6)--(2,4)--cycle; \path(1.5,4)--node{$i~i~i$}(1.5,6); 
				       \draw[<-](2.5,1)--(2.5,2.3); \path(2.6,2.3)--node{$m_{i,j+1}$}(2.6,4.7);  \draw[->](2.5,4.7)--(2.5,6);
\end{tikzpicture}
}}\cr
\end{array}
\end{equation} 
On the left are given the various constraints that may apply to entries in the $i$th row of $G$ for various $j$. 
These govern the entries $i$ that appear in the $j$th and $(j+1)$th rows of $S$ as illustrated on the right, where
each isolated $i$ in a box must appear, while the triples $i\,i\,i$ are intended to indicate optional sequences
of $i$s of various possible lengths.

Case (L) corresponds to the left-saturation of the betweenness condition (\ref{eqn-between}), and in this case there
are no entries $i$ in the $j$th rows of $S$ and thus no contribution to $\wgt(G)$. 
Case (R) corresponds to the right-saturation of (\ref{eqn-between}) and implies that there is at least one
entry $i$ in row $j$ of $S$ and this entry lies immediately above an entry $i$ in row $j+1$
as a result of the strictness condition (\ref{eqn-strict}) applied to entries $i$ in row
$j+1$. It follows that its contribution to $\wgt(S)$ is $(y_i-a_k)$ where $k=m_{i-1,j}$ is the number of steps 
it is from the main diagonal. This accounts for the term $\chi(R)(y_i-a_k)$ in (\ref{eqn-wgtG}).
The case (B) is the one, sometimes called special~\cite{Okada}, in which the betweennness condition is strict on both sides.
In this case there is at least one entry $i$ in row $j$ of $S$, but the leftmost such $i$ has no entry $i$ either to its left or 
vertically beneath it. Its contribution to $\wgt(S)$ is therefore $x_i+y_i$. 
This accounts for the term $\chi(B)(x_i+y_i)$ in (\ref{eqn-wgtG}).

As can be seen from the above diagrams, in cases (R) and (B) there may remain additional entries $i$ in row $j$ of $S$
and in both cases these are to the right of the leftmost $i$ that we have previously identified, and 
contribute to $\wgt(S)$ a contribution $(x_k+a_k)$ with $k$ equal to the number of steps to the right of the main diagonal. 
This is the origin of the terms $(x_i+a_k)$ in the
second line of (\ref{eqn-wgtG}). Finally the remaining factors of $(x_i+a_k)$ in (\ref{eqn-wgtG}) arise from the 
entries $m_{ii}$ in $G$ that specify a sequence of $m_{ii}$ entries $i$ in row $i$ of $S$ that start on the main diagonal,
with $k$ varying from $0$ to $m_{ii}-1$, as required to ensure that $\wgt(G)=\wgt(S)$, as required. 
\qed
\bigskip

Finally, we make contact with the results of Tokuyama~\cite{Tokuyama} and Bump et al.~\cite{BMN} that motivated
this work in the first place.

\begin{Corollary}~[Tokuyama]~\cite{Tokuyama}
\label{cor-TokG}
Let $\lambda=\mu+\rho$ with $\ell(\mu)\leq n$ and $\rho=(n-1,\ldots,1,0)$, then for $\x=(x_1,x_2,\ldots,x_n)$ and any $t$
\begin{equation}\label{eqn-TokG}
  \sum_{G\in{\cal G}^\lambda} t^{\#R(G)} (1+t)^{\#B(G)} \prod_{i=1}^n x_i^{\sum_{j=1}^{i} m_{ij}- \sum_{j=1}^{i-1} m_{i-1,j}}  
	=  \prod_{1\leq i<j\leq n} (x_i+tx_j)\  s_\mu(\x) \,,
\end{equation}
where $\#R(G)$ and $\#B(G)$ are the numbers of triples $(m_{ij},m_{i-1,j},m_{i,j+1})$ in $G$ satisfying the conditions
(R) and (B) of (\ref{eqn-LRB}), that is to say the number that are right-saturated and the number that are neither right nor left saturated,
respectively, and the exponent of $x_i$ is the difference between the sum of entries in the $i$th row of $G$ and the sum
of entries in the $(i-1)$th row of $G$, with the $0$th row defined to be empty.  
\end{Corollary}
\bigskip

\noindent{\bf Proof}:~~
This result, which makes precise Tokuyama's identity (\ref{eqn-Tok}), is a special case of Corollary~\ref{cor-G}. 
First, it should be noted that the difference between the 
use of $\lambda=\mu+\delta$ and $\lambda=\mu+\rho$ in Corollaries~\ref{cor-G} and \ref{cor-TokG}, respectively, just amounts to 
dropping the contribution $x_1x_2\ldots x_n$ that comes from the diagonal entries of $S$ and amounts to subtracting
$(1,1,\ldots,1)$ from $\lambda$. Then, the left hand side of (\ref{eqn-TokG}) is an immediate consequence of setting
$\a=(0,0,\ldots)$ and $y_i=tx_i$ for $i=1,2,\dots,n$ on the right hand side of (\ref{eqn-wgtG}) and collecting up the
terms in $t$, $(1+t)$ and $x_i$. Applying the same conditions to the right hand side of (\ref{eqn-tok-fact-G}) without the 
factor $x_1x_2\ldots x_n$ then yields the right hand side of (\ref{eqn-TokG}), as required.
\qed
\bigskip

Before proceeding to the next corollary it is convenient to introduce a small lemma
\begin{Lemma}
\label{lem-NE-SW}
Let $\lambda=\mu+\delta$ with $\mu$ a partition of length $\ell(\mu)\leq n$ and $\delta=(n,n-1,\ldots,1)$
and let $m=\lambda_1$. For $A\in{\cal A}^\lambda$ let $C$ be the corresponding compass point matrix and
let $\#\XY$ be the number of CPM entries $\XY$ in $C$. Then
\begin{equation}
    \#\SW = \#\NE +|\mu|\,.
\end{equation}  
\end{Lemma}

\noindent{\bf Proof}:~~
Let $\#\XY_i$ be the number of entries $\XY$ in the $i$th row of $C$
and let $\#i$ be the number of entries $i$ in the corresponding shifted tableau $S$. 
With this notation, 
\begin{equation}\label{eqn-lem1}
\#\NS_i+\#\NW_i+\#\NE_i = \sum_{k=1}^{i-1} \sum_{j=1}^m a_{kj} = i-1\,.
\end{equation}
Here the first step follows from the fact that the tabulation of (\ref{eqn-ASM-CPM-SIC}) implies
that the column sum in $A$ above the position of each entry $\XY$ in the $i$th row of $C$ 
is $1$ or $0$ according as $\XY$ is or is not in $\{\NS,\NW,\NE\}$, 
and the second step from the fact that the sum of entries in each row of $A$ is $1$.
However
\begin{equation}\label{eqn-lem2}
\#\WE_i+\#\NW_i+\#\SW_i = \#i\in S    \quad\hbox{and}\quad \#\WE_i=\#\NS_i+1
\end{equation}
since each entry $\XY\in\{\WE,\NW\SW\}$ in the $i$th row of $C$ gives rise to an entry $i$ in $S$
and each entry $\WE$ or $\NS$ in the $i$th row of $C$ corresponds to an entry $1$ or $\ov{1}$, respectively, 
in the $i$th row of $A$ whose row sum is $1$.
Combining these identities and summing over $i$ gives
\begin{equation}\label{eqn-lem3}
  \#\SW-\#\NE = |\lambda| - \sum_{i=1}^n i = |\mu|\,,
\end{equation}
as required.
\qed
\bigskip

This identity allows us to prove the following as a direct consequence of Corollary~\ref{cor-A}.

\begin{Corollary}~[Bump, McNamara and Nakasuji]~\cite{BMN}. 
\label{cor-TokBMN}
Let $\mu$ be any partition of length $\ell(\mu)\leq n$ and let $m=\mu_1+n$, then for $\z=(z_1,z_2,\ldots,z_n)$,
$\alpha=(\alpha_1,\alpha_2,\ldots )$ and any $t$, the partition function of the $6$-vertex planar spin configuration model 
takes the form
\begin{equation}\label{eqn-Zbeta}
   Z(\mathfrak{S}_{\mu,t}^{\Gamma}) = \sum_{\mathfrak{s}\in\mathfrak{S}_{\mu,t}^\Gamma} \prod_{i=1}^n\prod_{j=1}^m \beta_{ij}\,,  
\end{equation}
where the sum is over all possible internal spin states $\mathfrak{s}$ consistent with a given set of external spin states.
The six types of vertex at $(i,j)$ carry the Boltzmann weights $\beta_{ij}$ as tabulated below:
\begin{equation}\label{eqn-wgtC-Z}
\begin{array}{|c|c|c|c|c|c|c|c|}
\hline
\begin{array}{cc}\hbox{Spin states}\cr\hbox{at $(i,j)$}\end{array}
&
\vcenter{\hbox{
\begin{tikzpicture}[x={(0in,-0.21in)},y={(0.21in,0in)}] 
\draw(1,2)node{$\bullet$};
\draw[very thick,fill=black] (0.3,2) circle (0.04in); 
\draw[very thick,fill=white] (1.7,2) circle (0.04in); 
\draw[very thick,fill=white] (1,1.3) circle (0.04in); 
\draw[very thick,fill=black] (1,2.7) circle (0.04in); 
\draw[very thick] (1,2)to(0.47,2); 
\draw[very thick] (1,2)to(1.53,2); 
\draw[very thick] (1,2)to(1,1.47); 
\draw[very thick] (1,2)to(1,2.53); 
\end{tikzpicture}
}}
&
\vcenter{\hbox{
\begin{tikzpicture}[x={(0in,-0.21in)},y={(0.21in,0in)}] 
\draw(1,2)node{$\bullet$};
\draw[very thick,fill=white] (0.3,2) circle (0.04in); 
\draw[very thick,fill=black] (1.7,2) circle (0.04in); 
\draw[very thick,fill=black] (1,1.3) circle (0.04in); 
\draw[very thick,fill=white] (1,2.7) circle (0.04in); 
\draw[very thick] (1,2)to(0.47,2); 
\draw[very thick] (1,2)to(1.53,2); 
\draw[very thick] (1,2)to(1,1.47); 
\draw[very thick] (1,2)to(1,2.53); 
\end{tikzpicture}
}}
&
\vcenter{\hbox{
\begin{tikzpicture}[x={(0in,-0.21in)},y={(0.21in,0in)}] 
\draw(1,2)node{$\bullet$};
\draw[very thick,fill=white] (0.3,2) circle (0.04in); 
\draw[very thick,fill=white] (1.7,2) circle (0.04in); 
\draw[very thick,fill=black] (1,1.3) circle (0.04in); 
\draw[very thick,fill=black] (1,2.7) circle (0.04in); 
\draw[very thick] (1,2)to(0.47,2); 
\draw[very thick] (1,2)to(1.53,2); 
\draw[very thick] (1,2)to(1,1.47); 
\draw[very thick] (1,2)to(1,2.53); 
\end{tikzpicture}
}}
&
\vcenter{\hbox{
\begin{tikzpicture}[x={(0in,-0.21in)},y={(0.21in,0in)}] 
\draw(1,2)node{$\bullet$};
\draw[very thick,fill=black] (0.3,2) circle (0.04in); 
\draw[very thick,fill=black] (1.7,2) circle (0.04in); 
\draw[very thick,fill=black] (1,1.3) circle (0.04in); 
\draw[very thick,fill=black] (1,2.7) circle (0.04in); 
\draw[very thick] (1,2)to(0.47,2); 
\draw[very thick] (1,2)to(1.53,2); 
\draw[very thick] (1,2)to(1,1.47); 
\draw[very thick] (1,2)to(1,2.53); 
\end{tikzpicture}
}}
&
\vcenter{\hbox{
\begin{tikzpicture}[x={(0in,-0.21in)},y={(0.21in,0in)}] 
\draw(1,2)node{$\bullet$};
\draw[very thick,fill=white] (0.3,2) circle (0.04in); 
\draw[very thick,fill=white] (1.7,2) circle (0.04in); 
\draw[very thick,fill=white] (1,1.3) circle (0.04in); 
\draw[very thick,fill=white] (1,2.7) circle (0.04in); 
\draw[very thick] (1,2)to(0.47,2); 
\draw[very thick] (1,2)to(1.53,2); 
\draw[very thick] (1,2)to(1,1.47); 
\draw[very thick] (1,2)to(1,2.53); 
\end{tikzpicture}
}}
&
\vcenter{\hbox{
\begin{tikzpicture}[x={(0in,-0.21in)},y={(0.21in,0in)}] 
\draw(1,2)node{$\bullet$};
\draw[very thick,fill=black] (0.3,2) circle (0.04in); 
\draw[very thick,fill=black] (1.7,2) circle (0.04in); 
\draw[very thick,fill=white] (1,1.3) circle (0.04in); 
\draw[very thick,fill=white] (1,2.7) circle (0.04in); 
\draw[very thick] (1,2)to(0.47,2); 
\draw[very thick] (1,2)to(1.53,2); 
\draw[very thick] (1,2)to(1,1.47); 
\draw[very thick] (1,2)to(1,2.53); 
\end{tikzpicture}
}}
\cr
\hline
\beta_{ij}& 1 & (1+t)z_i & t &1 & z_i-t\alpha_j &z_i+\alpha_j\cr
\hline
c(i,j)&\WE &\NS & \NE & \SE & \NW & \SW \cr
\hline
\end{array}
\end{equation}	
where $\tikz \draw[very thick,fill=black] circle (0.04in);$ and $\tikz \draw[very thick,fill=white] circle (0.04in);$
signify spin up and down states, respectively. Then 
\begin{equation}\label{eqn-TokBMN}
    Z(\mathfrak{S}_{\mu,t}^{\Gamma})  =  \prod_{1\leq i<j\leq n} (tz_i+z_j)\  s_\mu(\z\,|\,\alpha) \,.
\end{equation}
\end{Corollary}
\bigskip

\noindent{\bf Proof}:~~
It should first be noticed that the $6$-vertex model spin state configurations are in bijective correspondence 
with SICs, ASMs and CPMs. The easiest way to implement the bijection between spin states $\mathfrak{s}$ and CPMs $C$ 
is to rotate $\mathfrak{s}$ through $\pi$ and map the resulting vertices to CPM entries $\XY$ as tabulated above.
It follows that 
\begin{equation}\label{eqn-Z-AC}
   Z(\mathfrak{S}_{\mu,t}^{\Gamma}) 
	  = \sum_{A\in{\cal A}^\lambda} \prod_{i=1}^n\prod_{j=1}^m \wgt(c_{ij}) 
\end{equation}
where $A$ is the ASM corresponding to the CPM $C$, $\lambda=\mu+\delta$ and $\wgt(c_{i,j})=\beta_{ij}$ for all $(i,j)$.	
This expression may then be evaluated				
by specialising the Boltzmann weights of (\ref{eqn-wgtC}) in such a way as to give those of (\ref{eqn-wgtC-Z}) modified
by moving the factor $t$ from $\NE$ to $\SW$ through the use of Lemma~\ref{lem-NE-SW}.
\begin{equation}\label{eqn-wgtC-Zmod}
\begin{array}{|c|c|c|c|c|c|c|c|c|}
\hline
&c_{ij}&\WE &\NS & \NE & \SE & \NW & \SW \cr
\hline
(\ref{eqn-wgtC-Z})~\hbox{modified}&\wgt(c_{ij})& 1 & (1+t)z_i & 1 &1 & z_i-t\alpha_j &t(z_i+\alpha_j)\cr
\hline
(\ref{eqn-wgtC})&\wgt(c_{ij})& 1 & x_i+y_i & 1 &1 & y_i-a_j &x_i+a_j\cr
\hline
\end{array}
\end{equation}							
Clearly these coincide if we set $x_i=tz_i$, $y_i=z_i$ and $a_i=t\alpha_i$ for $i=1,2,\ldots,n$. 
Comparing (\ref{eqn-wgtA}) and (\ref{eqn-tok-fact-A}) and remembering to include an overall factor 
$t^{-|\mu|}$ as required by Lemmma~\ref{lem-NE-SW}, we find
\begin{equation}
    Z(\mathfrak{S}_{\mu,t}^{\Gamma}) = t^{-|\mu|} \prod_{1\leq i<j\leq n} (tz_i+z_j)  \ s_\mu(t\z\,|\,t\alpha) 
\end{equation}
However the factorial Schur function $s_\mu(t\z\,|\,t\alpha)$ is homogeneous of degree $|\mu|$ in factors of the form $(tz_i+t\alpha_j)$,
so that  $t^{-|\mu|} s_\mu(t\z\,|\,t\alpha)=s_\mu(\z\,|\,\alpha)$, as required to complete the proof of (\ref{eqn-TokG}).
\qed
\bigskip

As a final corollary it is rather easy to recover the following result originally due to Lascoux~\cite{Lascoux} and rederived both by McNamara~\cite{McNamara} and 
by Bump {\it et al.}~\cite{BMN}:
\begin{Corollary}
Let $\x=(x_1,x_2,\ldots,x_n)$, $\a=(a_1,a_2,\ldots)$, $\delta=(n,n-1,\ldots,1)$, $\rho=(n-1,n-2,\ldots,0)$, $\mu$ be a partition 
of length $\ell(\mu)\leq n$, $\lambda=\mu+\delta$ and $\kappa=\mu+\rho$, with $\lambda'$ and $\kappa'$ the partitions conjugate to $\lambda$
and $\kappa$, respectively. Then writing $\z^\nu=z_1^{\nu_1}z_2^{\nu_2}\cdots$ for any $\z=(z_1,z_2,\ldots)$ and any $\nu=(\nu_1,\nu_2,\ldots)$, we have
\begin{equation}\label{eqn-Z-McN}
    Z(\mathfrak{S}_{\mu}) =  \frac{\x^\rho}{\a^{\kappa'}} (-1)^{|\kappa|} \ s_\mu(\x\,|\,\a)\,, 
\end{equation}
where $Z(\mathfrak{S}_\mu)$ is the partition function of the $6$-vertex model with Boltzman weights $\beta_{ij}=\wgt(c_{ij})$
given by
\begin{equation}\label{eqn-wgtC-McN}
\begin{array}{|c|c|c|c|c|c|c|c|}
\hline
c_{ij}&\WE &\NS & \NE & \SE & \NW & \SW \cr
\hline
\wgt(c_{ij})& 1 & -x_i/a_j & 1 &1 & 1 & -(x_i/a_j+1)\cr
\hline
\end{array}
\end{equation}							
\end{Corollary}

\noindent{\bf Proof}:~~
To see this one sets $y_i=0$ in the various $\wgt(c_{ij})$ taken from $(\ref{eqn-wgtC})$ of Corollary~\ref{cor-A}. This yields
\begin{equation}\label{eqn-AC-y0}
    \sum_{A\in{\cal A}^\lambda} \prod_{i=1}^n\prod_{j=1}^m \wgt(c_{ij}) = \x^\rho\ s_\mu(\x\,|\,\a)\,,
\end{equation}
with $\wgt(c_{ij})$ is given by the $y_i=0$ values specified below
\begin{equation}\label{eqn-wgtC-McNlinear}
\begin{array}{|c|c|c|c|c|c|c|c|c|}
\hline
&c_{ij}&\WE &\NS & \NE & \SE & \NW & \SW \cr
\hline
(\ref{eqn-wgtC})&\wgt(c_{ij})& 1 & x_i+y_i & 1 &1 & y_i-a_j &x_i+a_j\cr
\hline
y_i=0&\wgt(c_{ij})& 1 & x_i & 1 &1 & -a_j & x_i+a_j\cr
\hline
\end{array}
\end{equation}	
To effect the transition from (\ref{eqn-AC-y0}) to the required (\ref{eqn-Z-McN}) it is necessary to 
reassign the contributions $-a_j$ arising from each entry $NW$ in $C$.
This can be done by noting that if $\#\XY_j$ now represents the number of entries $\XY$ in the $j$th column of $C$
then $\#\WE_j+\#\NW_j+\#\SW_j$ is the number of entries in the $j$th diagonal of the 
corresponding unprimed shifted tableau $S$ of shape $\lambda$, but this number
of entries is $\lambda'_j$. It follows that
\begin{equation}
\begin{array}{rcl}
     (-a_j)^{\#\NW_j} &=& (-a_j)^{\lambda'_j-\#\WE_j-\#\SW_j} = (-a_j)^{\lambda'_j-\chi(j\in\{\lambda_1,\lambda_2,\ldots,\lambda_n\})-\#\NS_j-\#\SW_j}\cr
		    &=& (-a_j)^{\lambda'_{j+1}-\#\NS_j-\#\SW_j} = (-a_j)^{\kappa'_j-\#\NS_j-\#\SW_j} \cr
\end{array}				
\end{equation} 
where use has been made of the fact that $\#WE_j=\#NS_j+1$ or $\#NS_j$ according as
$j$ is or is not a part of $\lambda$, and the observation that for any strict partition $\lambda$ its
conjugate $\lambda'$ is such that $\lambda'_{j+1}=\lambda'_{j}-1$ or $\lambda'_{j}$ again according as 
$j$ is or is not a part of $\lambda$. This implies that we can pass from (\ref{eqn-AC-y0}) to
(\ref{eqn-Z-McN}) by changing the weights from the $y_i=0$ set in (\ref{eqn-wgtC-McNlinear}) to those
of (\ref{eqn-wgtC-McN}) and dividing on the right by the product over $j$ of $(-a_j)^{\kappa'_j}$, that is to say multiplying by $(-1)^{|\kappa|}/\a^{\kappa'}$.
\qed
\bigskip

\noindent{\bf Acknowledgements} \medskip

We have made extensive use of Maple to generate examples and verify results.
We are grateful to the organizers of the Banff International Research Station for Mathematical Innovation and Discovery (BIRS) 2010 5-day workshop, ``Whittaker Functions, Crystal Bases, and Quantum Groups'' where the first author (AMH) first learned of the problems in this area. She also acknowledges the support of a
Discovery Grant from the Natural Sciences and Engineering Research Council of
Canada (NSERC). The second author (RCK) is grateful for the hospitality and financial support
extended to him while visiting
Wilfrid Laurier University, Waterloo 
where some of this work was carried out. 
The authors thank Anna Bertiger for many helpful conversations.

\end{document}